\documentstyle[11pt]{article}
\date{}

\def\of{{\cal O}_{\varphi}}

\def\b{{\cal B}}
\def\e{{\cal E}}
\def\l{{\cal L}_{{\cal B}}(X)}

\def\sp{{\cal S}_p(X)}
\def\s1{{\cal S}_1 (X)}
\def\supp{\hbox{supp}}

\def\f{\varphi}

\def\ben{\begin{enumerate}}
\def\een{\end{enumerate}}

\def\zR{{\rm I}\!{\rm R}}

\def\zC{{\rm C}\!\!\!\vrule height 6.4pt depth -0.4pt width
1pt \ }
\def\zN{{{\rm I}\!{\rm N}}}

\def\diagram#1{{\normallineskip=2pt
\normalbaselineskip=0pt\matrix{#1}}}

\def\hrarr^#1_#2{ \mathrel{
\mathop{\hbox to .3in{\rightarrowfill}}
\limits^{\scriptstyle#1}_{\scriptstyle#2}  }}

\def\hlarr^#1_#2{ \mathrel{
\mathop{\hbox to .5in{\leftarrowfill}}
\limits^{\scriptstyle#1}_{\scriptstyle#2}  }}

\def\vdarr#1#2{\llap{$\scriptstyle#1$}
\left\downarrow\vcenter to .3in{}\right.
\rlap{$\scriptstyle#2$}}

\def\ddrarr#1{\searrow\displaystyle\rlap{$\vcenter
{\hbox{$\scriptstyle #1$}}$}}

\newtheorem{lemm}{Lemma}[section]
\newtheorem{theo}[lemm]{Theorem}
\newtheorem{prop}[lemm]{Proposition}
\newtheorem{coro}[lemm]{Corollary} 
\newtheorem{rema}[lemm]{Remark} 
 
\newtheorem{exam}[lemm]{Example}

\newtheorem{hipo}[lemm]{Hypothesis}

\newenvironment{proof}[1]{
  \trivlist \item[\hskip \labelsep{\bf #1}]}{\hfill\mbox{$\Box$}
  \endtrivlist}

\begin{document}
%\fontsize {14}{19}\selectfont
%\Large

\title{ HOMOTOPY OF VECTOR STATES}
%\centerline{\bf \Large VECTOR STATES IN C*-MODULES\footnote{1991 
%				Mathematics Subject Classification: 46L30, 46L05, 46L10.}}
\vskip.6cm

\author{ Esteban Andruchow and Alejandro Varela}
\vskip.7cm  
\maketitle
\centerline{\large Instituto de Ciencias -- Univ. Nac. de Gral. Sarmiento}
\centerline{\large Argentina}
\vskip1cm

%\centerline{\today}

\vskip0.5cm
\abstract{ Let $\b$ be a C$^*$-algebra and $X$ a C$^*$ Hilbert 
$\b$-module. If  $p\in \b$ is a projection, denote by 
$\sp =\{x\in X : \langle x,x\rangle =p\}$, the $p$-sphere of $X$. 
For $\f$ a state of $\b$ with support $p$  in $\b$ and $x\in \sp$, 
consider the state $\f_x$ of $\l$ given by $\f_x(t)=\f(\langle x,t(x)\rangle )$. 
In this paper we study certain sets associated to these states, 
and examine their topologic properties. As an application of these techniques, we prove that the space of states of the hyperfinite II$_1$ factor ${\cal R}_0$, with support equivalent to a given projection $p\in {\cal R}_0$, regarded with the norm topology (of the conjugate space of ${\cal R}_0$), has trivial homotopy groups of all orders.

The same holds for the space
$$
{\cal S}_p({\cal R}_0)=\{v\in {\cal R}_0:v^*v=p\}\subset {\cal R}_0
$$
of partial isometries with initial space $p$, regarded with the ultraweak topology.
 }

\bigskip


\bigskip
%\fontsize {14}{16}\selectfont
\section{Introduction}
Tools from homotopy theory have been used in operator algebras for quite some time. Starting with N.H. Kuiper's theorem \cite{kui}, establishing the contractibility of the unitary group of an infinite dimensional Hilbert space, following with further generalizations, to properly infinite von Neumann algebras (\cite{bre}, \cite{ass}, \cite{bru}). The case of a finite von Neumann algebra was considered by Araki,  M. Smith and L. Smith, who  showed in \cite{ass} that the $\pi_1$ group of the unitary group of a II$_1$ factor is isomorphic to the additive group of the reals $\zR$. This result was later extended by Schr\"oder in \cite{sch}. Also some results appeared computing the homotopy type of the unitary groups of certain classes of C$^*$-algebras (\cite{han}, \cite{zha}).
The topology considered for the unitary groups of the von Neumann algebras in these papers is the one induced by the norm of the algebra. Only a few years ago Popa and Takesaki \cite{pota} studied the homotopy theory of the unitary and automorphism groups of a factor in the weak topologies. 
 
In this paper we focus on the homotopy type of certain sets of states, namely the states having support equivalent to a fixed projection, both with uniform (norm) and weak topologies.

Let $\b$ be a von Neumann algebra. The set of faithful and normal states of a $\b$ is a convex set, and therefore is contractible. If the algebra is finite, this set is the same as the set of states with support {\it equivalent} to the identity $1\in \b$. This paper considers the question of what happens when one replaces $1$ by a proper projection $p\in \b$. In other words, what is the homotopy type of the set of normal states with support (Murray-von Neumann) equivalent to $p$.

It turns out that this set is closely related to the set of partial isometries of $\b$ with initial space $p$, i.e. the elements $v\in \b$ such that $v^*v=p$. Namely, by means of the map
$$
v \times \f \mapsto \f(v^* \cdot v)
$$
where $\f$ is a given state with support $p$.
The purpose of this paper is the study of the properties of this map: mainly, under which assumptions it is a fibration. And in the affirmative case, to use this fibration to compute the homotopy type of the set of states described.

First, there is the question of what is the right topology to consider in the set of partial isometries (among the various topologies available in $\b$). It turns out the norm topology
of $\b$ forces on the set of states, via the map above, a topology stronger than the norm topology for functionals. In order to have fibration properties for this map, and have it induce the usual norm topology on the set of states, one has to consider on (the set of isometries of) $\b$ the ultraweak topology. The paper by S. Popa and M. Takesaki \cite{pota}, deals with the topologic properties of the unitary group in the weak topologies. There Michael's theory of continuous selection \cite{mic} is used in a remarkable manner, to obtain cross sections (i.e. fibration properties) for the quotient of the unitary groups of an inclusion of factors. We shall profit from their technique in our context.

In sections 3 and 4 of this paper we shall put our problem in the broader context of Hilbert C$^*$-modules. All the concepts above have their analogue in this general setting. We shall therefore state first the results valid in the general context. In sections 5 and 6 we shall return again to our original situation, and shall be able to state the result described in the abstract, for a class of II$_1$ factors considered in \cite{pota}, which encloses the hyperfinite factor ${\cal R}_0$.

Let now $\b$ be a C$^*$- algebra, $X$ a right C$^*$-module over  $\b$, and $\l$ the C$^*$-algebra of adjointable operators of $X$. If $p\in \b$ is a projection, denote by $\sp=\{x\in X: \langle x,x\rangle =p \}$ 
the $p$-sphere of $X$. We shall study the states of $\l$ which are {\it pure} in the modular sense. That is, for a state $\f$ of $\b$ and a vector $x\in \sp$, we consider the state
$\f_x$ with density $x$, given by
$$
\f_x(t)=\f(\langle x,t(x)\rangle) \ \ t\in \l .
$$
We require that the state $\f$ has support projection $p$. 

If $x,y \in X$, let
$\theta_{x,y} \in \l$ be the ``rank one'' operator given by 
$\theta_{x,y}(z)=x\langle y,z\rangle $. If $\langle x,x\rangle =p$ then the operator $\theta_{x,x}=e_x$ is a selfadjoint projection, and all projections arising in this manner, from vectors on $\sp$, are mutually (Murray-von Neumann) equivalent. It turns out that these modular, or {\it vector states} as we shall subsequently call them, are precisely the states of $\l$ with support of rank one, i.e. equal to one of these projections $e_x$.

Namely we  are interested in the following sets of states,  
$$
\of =\{ \f_x : x\in \sp\}
$$
for $\f$ a fixed state in $\b$, with support projection $\supp(\f)=p$ in $\b$,
$$
\Sigma_p(\b)=\{ \hbox{states of } \b \hbox{ with support } p\}
$$
and
$$
\Sigma_{p,X}=\{\psi_x : \psi \in \Sigma_p(\b) , x \in \sp \}.
$$
The (convex) set $\Sigma_p(\b)$ is considered with the relative topology induced by the usual norm of the conjugate space of $\b$. The other two sets $\of$ and $\Sigma_{p,X}$ are endowed  with two types of natural topologies. 

First, in sections 3 and 4, with norm related topologies, given by

\begin{enumerate}
\item{-} $\of$ with the metric $d_\f$, $d_\f(\Phi, \Psi)=\inf \{\|x-y\| : \f_x=\Phi, \f_y=\Psi \}$. 
\item{-} $\Sigma_{p,X}$ with the metric $d$, $d(\Phi,\Psi)=\|\Phi -\Psi\| + 
		\|\supp(\Phi)-\supp(\Psi)\|$.
\end{enumerate}
These metrics $d$ and $d_\f$ do come up naturally if one looks for continuity of the maps
$$
\sp \to \of \ , \ x\mapsto \f_x
$$
and 
$$
\sp \times \Sigma_p(\b) \to \Sigma_{p,X} \ , \ (x,\f) \mapsto \f_x .
$$
It turns out that these maps are also fibrations,
a fact which enables one to
examine the homotopy theory  of these spaces, comparing it with the homotopy of the spheres $\sp$ studied in \cite{acs} and \cite{acs2}
For example, we obtain 
that if $p\b p$ is a von Neumann algebra and $Xp$ is  selfdual (as a 
$p\b p$-module), then $\of$ is simply connected. Or that if $p\b p$ 
is a properly infinite von Neumann algebra, and again $Xp$ is selfdual, 
then $\Sigma_{p,X}$ has trivial homotopy groups of all orders. 

Second, in section 5, the sets $\of$ and $\Sigma_{p,X}$ are considered with weaker topologies, related to the w$^*$ topology of the module $X$, (here $\b$ is supposed to be a von Neumann algebra, and $X$ a selfdual module \cite{pas}).
This is done by means of the natural representation $\rho$ of $\l$ in the Hilbert space ${\cal H}$ based on the algebraic tensor product $X\otimes H$. We choose $H$ the Hilbert space of a standard representation $\b \hookrightarrow B(H)$. This representation  $\rho$ was introduced in \cite{rie}, and studied also in \cite{pas}. Here we show that via $\rho$, the states in the set $\Sigma_{1,X}$ are pure in the classic sense, i.e. they are induced  by vectors of the form $x\otimes \xi$ of ${\cal H}$. We denote by ${\cal A}(X)$ the set of all such density vectors. $\Sigma_{1,X}$ appears then as the set $\Omega_{{\cal A}(X)}=\{\omega_{x\otimes \xi}: x\otimes \xi \in {\cal A}(X)\}$, where $\omega_{x\otimes \xi}$ is the state induced by $x\otimes \xi$. Moreover $\s1 \times \Sigma_1(\b)$ is in one to one correspondence with ${\cal A}(X)$, $(x,\f)\leftrightarrow x\otimes \xi$ ($\f=\omega_\xi$), and the map $(x,\f)\mapsto \f_x$ translates as
$$
{\cal A}(X) \to \Omega_{{\cal A}(X)} , \ x\otimes \xi \mapsto \omega_{x\otimes \xi} .
$$
We endow ${\cal A}(X)\subset {\cal H}$ with the Hilbert space norm topology and $\Omega_{{\cal A}(X)}$ with the quotient topology induced by this map, which turns out to be the relative norm topology of the conjugate space of $\l$. We deal with $\b={\cal R}$ a II$_1$ factor with separable predual and the property that the tensor product ${\cal R}\otimes B(H)$ has a one parameter group of automorphisms which scales the trace \cite{pota} (a class which contains 
the hyperfinite II$_1$ factor). It is shown that the map $x\otimes \xi \mapsto \omega_{x\otimes \xi}$ is a trivial principal bundle (with structure group the unitary group $U_{\cal R}$ in the strong operator topology). This implies that in this case the homotopy groups of $\Omega_{{\cal A}(X)}$ coincide with the groups of the sphere $\s1$ in the (relative ) w$^*$ topology of $X$ (which is assumed selfdual). 

Using these results, in section 6,
we obtain that the set of states of such ${\cal R}$ with support equivalent to a fixed $p$ has trivial homotopy groups of all orders. Also it is proven that the set ${\cal S}_p({\cal R})$ of partial isometries of ${\cal R}$ with initial space $p$, regarded with the ultraweak topology, has trivial homotopy groups.

Another application of these results enables one to relate the homotopy groups of the unitary orbit ${\cal U}_\f=\{\f\circ Ad(u): u\in U_{{\cal R}}\}$ of $\f$ in the norm topology, with the homotopy groups of 
the unitary group of the centralizer ${\cal R}^\f$ in the strong operator topology. In particular one obtains that ${\cal U}_\f$ is simply connected.

\section{Preliminaries and notations}

Let us establish some basic facts and notations about the vector 
states $\f_x$.

We shall be concerned with states of $\b$ that have
their support in $\b$, a fact which holds automatically if $\b$ is a von Neumann algebra. 
\medskip

Each element $x\in \sp$ gives rise to a (non unital) $\ast$-isomorphism 
$i_x:p \b p \to  \l$,   $i_x (a)=\theta_{xa,x}$. 

Fix $x_0 \in \sp$. Let us recall from \cite{acs2} the principal fibre bundle, 
which we shall call the projective bundle
$$
\rho : \sp \to \e_{e_{x_0}}=\{ \hbox{projections in $\l$ 
equivalent to }e_{x_0}\}
$$
given by $\rho (x)=e_x$. Note that $\e_{e_{x_0}}$ depends only on $p$ and
 not on the choice of $x_0$ (all projections on $\e_{e_{x_0}}$ are of the 
 form $e_x$ for some $x\in \sp$). The structure group of the projective 
 bundle is the unitary group $U_{p\b p}$ of $p\b p$.

As is usual notation, if $\f$ is a faithful state of $\b$, $\b^\f$ is the
 centralizer algebra of $\f$, i.e. 
 $\b^\f=\{a\in \b: \f(ab)=\f(ba)\  \hbox{for all } b\in \b\}$. 
 If the support $\supp(\f)=p <  1$, then denote by $\b_p^\f$ the centralizer
  of the restriction of $\f$ to the reduced algebra $p\b p$.

Typically $a,b,c$ will denote elements of $\b$, $x,y,z$ elements of $X$ 
and $r,s,t$ elements 
of $\l$. $\b''$ will denote the von Neumann enveloping algebra of $\b$, 
and $X'$ the selfdual completion of $X$, which is 
a C$^*$-module over $\b''$ (\cite{pas}). By fibre bundle we mean a locally trivial fibre bundle, and by fibration we mean a surjective map having the homotopy lifting property (\cite{Sw}).  

\begin{lemm} 
Let $\f$ be a state of $\b$ with $\supp(\f)=p \in \b$, and $x$ an element in $\sp$. 
Then $\supp(\f_x)=e_x$.
\end{lemm}
\begin{proof}{Proof.}
Clearly $\f_x(e_x)=\f(\langle x,e_x(x)\rangle )=\f(\langle x,xp\rangle )=\f(p)=1$. Put $r=\supp(\f_x)\in {\cal L}_{\b''}(X')$. 
We have $r\le e_x$.  In particular $e_xr=re_x=r$. This implies that $r$ is of the 
form $e_y=\theta_{y,y}$, namely,
$y=r(x)\in X'$. Now $\langle r(x),r(x)\rangle =q$ is a projection in $\b''$, with $q\le p$. 
Indeed, 
$$
\begin{array}{rl}
\langle r(x),r(x)\rangle \langle r(x),r(x)\rangle & =
\langle r(x),r(x)\langle r(x),r(x)\rangle \rangle \\&\\
& = \langle r(x),\theta_{r(x),r(x)}(r(x))\rangle =\langle r(x),r(x)\rangle .
\end{array}$$
And $\langle r(x),r(x)\rangle p=\langle r(x),r(xp)\rangle =\langle r(x),r(x)\rangle $, i.e. $q\le p$. 
Now it is clear that $\f(q)=\f(\langle r(x),r(x)\rangle )=\f(\langle x,r(x)\rangle )=\f_x(r)=1$, 
which  implies that $q=p$, Therefore 
$\langle r(x)-x,r(x)-x\rangle =\langle r(x),r(x)\rangle +\langle x,x\rangle -\langle r(x),x\rangle -\langle x,r(x)\rangle =0$, 
since all these products equal $p$ 
(because $\langle r(x),x\rangle =\langle r^2(x),x\rangle =\langle r(x),r(x)\rangle $). 
Finally, $r(x)=x$ implies that $r=e_{r(x)}=e_x$.
\end{proof}
 
\begin{lemm}\label{soporte}
Let $\Phi$ be a state of $\l$ with $\supp(\Phi)=e_x$ for some $x\in \sp$. 
Then $\Phi=\f_x$ for  $\f$ a state in $\b$ with $\supp(\f)=p$. 
Namely $\f(a)=\Phi(i_x(a))$.
\end{lemm}
\begin{proof}{Proof.}
Put $\f=\Phi\circ i_x$ as above. First note that if $t\in \l$, then 
$e_x t e_x=\theta_{x\langle x,t(x)\rangle ,x}$.  
Then $\f_x(t)=\Phi(i_x(\langle x,t(x)\rangle ))=\Phi(\theta_{x\langle x,t(x)\rangle ,x})=
\Phi(e_x t e_x)=\Phi(t)$. It remains to see that $\supp(\f)=p$. 
Clearly $\f (p)=\Phi (\theta_{xp,x})=\Phi(e_x)=1$. Suppose that $q \le p$ is a 
projection in $\b''$ with $\f(q)=\Phi(\theta_{xq,x})=1$ ($\f$ here denotes the 
normal extension of the former $\f$ to $\b''$). 
Note that $\theta_{xq,x}=\theta_{xq,xq}=e_{xq}$ is in fact a projection 
(associated to $xq \in {\cal S}_q(X')$), and verifies $e_{xq}\le e_x$. 
It follows that $\theta_{xq,xq}=\theta_{x,x}$. 
Then $xq=\theta_{xq,xq}(x)=\theta_{x,x}(x)=x$, and therefore $q=p$.
\end{proof}
\begin{rema}\label{normal} \rm
If $\b$ is a von Neumann algebra,
the inner product of $X$ is weakly continuous and the state $\Phi$ 
of the preceding result is normal, then $\f=\Phi \circ i_x$ is also normal.
\end{rema}

%Suppose now that $\f$ is a state of $\b$ with 
%$\supp(\f)=p\in\b$, and $x,y \in \sp$. 
%The next result establishes when $\f_x=\f_y$.

\begin{prop}\label{fibra de fix}
Let $\psi, \f\in \Sigma_p(\b)$, $x,y\in \sp$. Then 
\ben

\item[a) ] $\f_x=\psi_x$ if and only if $\f=\psi$.
\item[b) ] $\f_x=\psi_y$ if and only if $\psi=\f\circ Ad(u)$, with $y=xu$ and $u\in U_{p\b p}$.
\item[c) ] $\f_x=\f_y$ if and only if $y=xv$, 
for $v$ a unitary element in $\b_p^\f$.
\een
\end{prop}

\begin{proof}{Proof.}
Let us start with a):  
$\f(b)=\f_x(\theta_{xb,x})=\psi_x(\theta_{xb,x})=\psi(b)$.

\smallskip
To prove b), suppose that $\f_x=\psi_y$.
Then they have the same support, i.e. $e_x=e_y$, which implies that there 
exists a unitary element $u\in U_{p\b p}$ such that $y=xu$ (see \cite{acs2}). Then
$$
\f_x(t)=\psi_y (t)= \psi(\langle xu,t(xu)\rangle )=\psi(u^*\langle x,t(x)\rangle u)=[\psi \circ Ad(u^*)]_x(t).
$$ 
Using part a), this implies that $\f=\psi \circ Ad(u^*)$, or 
$\psi=\f \circ Ad(u)$. 

\smallskip
To prove c), use b), and note that the unitary element $u \in U_{p\b p}$ satisfies $\f=\f \circ Ad(u)$, i.e. $u\in \b_p^\f$.
\end{proof}

\section{The  set $\of$}
In this section we shall consider the set $\of=\{ \f_x : x\in \sp\} $ 
for a fixed state $\f$ of $\b$ with $\supp (\f)=p$. Note that in the 
particular case when $X=\b$ is a finite von Neumann algebra and $p=1$,
then $\of$ is just the unitary orbit of $\f$. In the general case, 
there is a canonical map 
$$
\sigma : \sp \to \of , \ \sigma (x)=\f_x.
$$
Let us consider the following natural metric in $\of$:
$$
d_\f(\f_x, \f_y)= \inf \{ \|x'-y'\| : x',y' \in \sp , \f_{x'}=\f_x , \f_{y'}=\f_y \}
$$
It is clear that this metric induces the same topology as the quotient topology 
given by the map $\sigma$, also, that in view of \ref{fibra de fix} 
it can be computed as follows:
$$
d_\f(\f_x, \f_y)= \inf \{ \|x-yv\|: v \hbox{ unitary in } \b_p^\f\}.
$$

First note that this is indeed a metric. For instance, if $d_\f(\f_x, \f_y)=0$ 
then there exist unitaries $v_n$ in $\b_p^\f$ such that $\|x-yv_n\| \to 0$, 
i.e. $yv_n \to x$ in $\sp$. In particular $yv_n$  is a Cauchy sequence, and 
therefore $v_n$ is a Cauchy sequence, converging to a unitary $v$  in $\b_p^\f$. 
Then $x=yv$ and $\f_x=\f_y$. The other properties follow similarly.

With this metric, $\of$ is homeomorphic to the quotient $\sp / U_{\b_p^\f}$. 
The following result implies that the inclusion $\of \subset \b^*$ 
(=conjugate space of $\b$) is continuous.
\begin{lemm}
If $x,y \in \sp$, then $\|\f_x - \f_y\| \le 2\| x-y\|$. In particular
$$
\|\f_x - \f_y\| \le 2d_\f(\f_x, \f_y)
$$
where the norm $\| \ \|$ of the functionals denotes the usual norm of the conjugate space $\b^*$
\end{lemm}
\begin{proof}{Proof.}
If $t\in \l$, then 
$|\f_x(t)-\f_y(t)|\le |\f(\langle x,t(x-y)\rangle |+|\f(\langle x-y,ty\rangle )|$. Now by the Cauchy-Schwarz 
inequality $\|\langle x,t(x-y)\rangle \| \le \|t\| \ \|x-y\|$, and $\|\langle x-y,ty\rangle \| \le \|x-y\| \ \|t\|$. Then
$\|\f_x(t)-\f_y(t)\| \le 2\|t\| \ \|x-y\|$, and the result follows.
\end{proof}
We want  the map $\sigma :\sp \to \of$ to be a locally trivial 
fibre bundle. In order to obtain it we make the following assumption:

\begin{hipo}\label{hipotesis}
There exists a conditional expectation $E_\f : p\b p \to \b_p^\f$.
\end{hipo}
This is the case if for example $\b$ is a von Neumann algebra and $\f$ is normal.
 For the remaining of this section, we suppose that \ref{hipotesis} holds.
\begin{theo}
The map $\sigma :\sp \to \of$, $\sigma (x)=\f_x$ is a locally trivial fibre bundle. 
The fibre of this bundle is the unitary group $U_{\b_p^\f}$ of $\b_p^\f$.
\end{theo}
\begin{proof}{Proof.}
Since the spaces are homogeneous spaces,
it suffices to show that there exist continuous local cross sections at every 
point $x_0$ of $\sp$. Suppose that $d_\f(\f_x,\f_{x_0})< r< 1$. 
Then there exists a unitary operator $v$ in $\b_p^\f$ such that $\|xv - x_0\|< 1$. Then
$$
\|p - \langle xv,x_0\rangle \|=\|\langle x_0,x_0\rangle -\langle xv,x_0\rangle \|=\|\langle x_0-xv,x_0\rangle \|\le \|x_0 -xv\|< 1
$$
It follows that $\langle xv,x_0\rangle $ is invertible in $p\b p$.
 Therefore, one can find $r$ such that also 
 $E_\f(\langle xv,x_0\rangle )=v^*E_\f(\langle x,x_0\rangle )$ is invertible in $\b_p^\f$. 
 Then $E_\f(\langle x,x_0\rangle )$ is invertible. Let us put
$$
\eta_{x_0} (\f_x) = x \mu (E_\f(\langle x,x_0\rangle ))
$$
defined on $\{\f_x : d_\f(\f_x,\f_{x_0})< r\}$, where $\mu$ denotes the
 unitary part in the polar decomposition (of invertible elements) in 
 $\b_p^\f$: $c=\mu (c) (c^*c)^{1/2}$. 
First note that $\eta_{x_0}$ is well defined. If $x'$ is a vector in the 
fibre of $\f_x$, then $x'=xv$ for $v\in U_{\b_p^\f}$. Then 
$x' \mu(E_\f(\langle x',x_0\rangle ))=xv \mu (v^*E_\f (\langle x,x_0\rangle ))=x \mu (E_\f (\langle x,x_0\rangle ))$, 
where the last equality holds because $\mu (ua)=u\mu (a)$ if $u$ is unitary.
Next, $\eta_{x_0}(\f_{x_0})=x_0 \mu (E_\f(\langle x_0,x_0\rangle ))=x_0$, and $\eta_{x_0}$ 
is a cross section for $\sigma$, because $\mu (E_\f(\langle x,x_0\rangle ))$ is a unitary 
in $\b_p^\f$. Finally, let us see that $\eta_{x_0}$ is continuous. Suppose 
that $d_\f(\f_{x_n}, \f_x) \to 0$, then there exist unitaries $v_n$ in $\b_p^\f$ 
such that $x_n v_n \to x$. Then by the continuity of the operations,
$x_n v_n \mu(E_\f(\langle x_n v_n, x_0\rangle ))=x_n \mu (E_\f (\langle x_n,x_0\rangle )) \to x \mu (E_\f(\langle x,x_0\rangle ))$, 
i.e., $\eta_{x_0}(x_n) \to \eta_{x_0}(x)$.
It is clear from \ref{fibra de fix} that the fibre is $U_{\b_p^\f}$. Namely, $\sigma^{-1}(\f_x)=\{xv: v\in U_{\b_p^\f}\}$. Note that $xv_n \to xv$ in $\sigma^{-1}(\f_x)\subset \sp$ if and only if $v_n\to v$ in $U_{\b_p^\f}$.
\end{proof}
We shall use the following result, which is a straightforward fact from the theory of fibrations
\begin{lemm} \label{lemita}
Suppose that one has the following commutative diagram
$$
\diagram{
E &  \hrarr^{\pi_1}_{} & X  \cr
{} & \ddrarr{\pi_2 }{}       & \vdarr{}{p}\cr
{}  &      {}        &  Y ,
}
$$
where $E, X, Y$ are topological spaces, $\pi_1$, $\pi_2$ 
are fibrations and $p$ is continuo\-us and surjective. Then $p$ is also a fibration.
\end{lemm}

There is another natural bundle associated to $\of$, which is 
the mapping
$$
\of \to \e_e , \ \f_x \mapsto e_x ,
$$
where $e$ is any projection of the form $e_{x_0}$ for some $x_0 \in \sp$. 
Since $e_x=\supp (\f_x)$, we shall call this map $\supp$. In general, 
taking support of positive functionals does not define a continuous map. 
However it is continuous in this context, i.e. restricted to the set $\of$ 
with the metric $d_\f$. Indeed, as seen before, convergence of 
$\f_{x_n} \to \f_x$ in this metric implies the existence of unitaries 
$v_n$ of $\b_p^\f \subset p\b p$ such that $x_n v_n \to x$ in $\sp$. 
This implies that 
$e_{x_n v_n}= e_{x_n} \to e_x$. Moreover, one has
\begin{theo}
The map $\supp : \of \to \e_e$ is a fibration with 
fibre $U_{p \b p}/ U_{\b_p^\f}$. One has the following commutative diagram of fibre bundles
$$
\diagram{
\sp &  \hrarr^{\rho}_{} &\of  \cr
{} & \ddrarr{\sigma }{}       & \vdarr{}{\supp}\cr
{}  &      {}        &  \e_e .
}
$$
\end{theo}
\begin{proof}{Proof.}
That the diagram commutes is apparent. 
Since $\rho$ and $\sigma$ are fibre bundles, it follows using \ref{lemita} that $\supp$ is a fibration. Note that if $e_x=e_y$, then there exists $u\in U_{p\b p}$ such that $y=xu$. Then $\f_y=\f_{xu}= (\f \circ Ad(u^*))_x$, therefore $\supp^{-1}(e_x)=\{(\f\circ Ad(u^*))_x: u\in U_{p\b p} \}$. In \ref{fibra de fix} it was shown that $\f_x=\psi_x$ (where $\f$ and $\psi$ are states of $\b$ with support $p$) implies $\f=\psi$. So $\{\f \circ Ad(u): u\in U_{p\b p}\}$ parametrizes the fibres of $\supp$, and clearly this set is in one to one correspondence with $U_{p\b p}/U_{\b^\f_p}$. Now $(\f \circ Ad(u_n))_x\to (\f\circ Ad(u) )_x$ in $\of$ if and only if $\inf_{v\in U_{\b^\f_p}}\|xu_n - xuv\|=\inf_{v\in U_{\b^\f_p}}\|u_n-uv\|$, i.e. the class of $u_n$ tends to the class of $u$ in $U_{p\b p}/U_{\b^\f_p}$ (with the quotient topology induced by the norm of $\b$).
\end{proof}

One can use the homotopy exact sequences of these bundles to relate the homotopy 
groups of $\of$, $\sp$, $\e_e$, $U_{p\b p}$, $U_{\b_p^\f}$ and $U_{p\b p} / U_{\b_p^\f}$. 
Namely:
$$
\dots \pi_n(U_{\b_p^\f},p)\to \pi_n(\sp,x_0)\stackrel{\sigma_*}{\to}  
\pi_n(\of,\f_{x_0})\to \pi_{n-1}(U_{\b_p^\f},p)\to \dots
$$
where $x_0$ is a fixed element in $\sp$, and
%\fontsize {13}{19}\selectfont
$$
\dots \pi_n(U_{p\b p}/U_{\b_p^\f},[p])\to \pi_n(\of,\f_{x_0})
\stackrel{\supp_*}{\to} \pi_n({\cal E},e_{x_0})\to \pi_{n-1}(U_{p\b p}/U_{\b_p^\f},[p])
\to \dots
$$
%\fontsize {14}{19}\selectfont
with $\f$ a fixed state in $\Sigma_p(\b)$.

The first result uses simply the fact that $\sigma$ is continuous and surjective.
\begin{coro}
If $p\b p$ is a finite von Neumann algebra, then $\of$ is arcwise connected.
\end{coro}
\begin{proof}{Proof.}
If $p \b p$ is finite, it was shown in $\cite{acs2}$ that $\sp$ is connected.
\end{proof}
\begin{coro}
If $p \b p$ is a von Neumann algebra and the restriction of $\f$  to $p \b p$ 
is normal, then                        
$$
\pi_1(\of, \f_x) \cong \pi_1 (\e_e, e_x).
$$
If moreover $Xp$ is selfdual, then $\pi_1(\of, \f_x)=0$.
\end{coro}
\begin{proof}{Proof.}
The proof follows by applying the tail of the homotopy exact sequence of 
the bundle $\supp$, recalling from $\cite{av2}$ that the fibre 
$U_{p\b p} / U_{\b_p^\f}$ is simply connected. In the selfdual case, 
it was proven in $\cite{acs2}$ that the connected components of $\e_e$ 
are simply connected.
\end{proof}
\begin{coro} \label{2.7}
If $p \b p$ is a von Neumann algebra, $\f$ restricted to $p \b p$ is 
normal and $Xp$ is selfdual, then, for a fixed $x_0 \in \sp$, the inclusion 
map 
$$
i : U_{\b_p^\f} \hookrightarrow \sp , \ v\mapsto x_0v
$$
induces an epimorphism
$$
i_*: \pi_1(U_{\b_p^\f},p) \to \pi_1 (\sp, x_0).
$$
\end{coro}
\begin{proof}{Proof.}
This time use the homotopy exact sequence of $\sigma$, and the fact that 
in this case $\pi_1(\of,\f_{x_0})=0$.
\end{proof}
In other words, this result says that regardless of the size of the selfdual 
module $X$, any closed continuous curve $x(t) \in \sp$ with $x(0)=x(1)=x_0$ is 
homotopic to a closed curve of the form $x_0(t)=x_0 v(t)$, with $v(t)$ a curve 
of unitaries in $\b_p^\f$ with $v(0)=v(1)=p$.
\begin{coro}
Suppose that $X$ is selfdual. If either 
\begin{enumerate}
\item[a) ]  $p\b p$ is a properly infinite von Neumann algebra, 

\noindent or 
\item[b) ] $p\b p$ is a von Neumann algebra of type II$_1$ with $\l$ properly infinite,
\end{enumerate}
then for $n\geq 1$
$$
\pi_n (\of,\f_{x_0})\cong \pi_{n-1}(U_{\b_p^\f},p).
$$
\end{coro}
\begin{proof}{Proof.}
In both  cases, a) and b), one has that $\sp$ is contractible (see \cite{acs2}). 
Therefore the proof follows writing down the homotopy exact sequence of the 
fibre bundle $\sigma$.
\end{proof}

Situation b) occurs for example if $p\b p$ is a II$_1$ factor and $Xp$ 
is not finitely generated over $p\b p$.

Finally let us state an analogous result for general C$^*$-algebras $\b$ 
(under the hypothesis \ref{hipotesis}) for the module 
$X=H_{\b}=\b \otimes \ell^2$. Here we use the fact  \cite{acs2}, 
that ${\cal S}_p(H_\b)$ is contractible. The proof follows similarly.
\begin{coro}
If \ref{hipotesis} holds, and $X=H_\b$, then for $n\ge 1$
$$
\pi_n (\of,\f_{x_0})\cong \pi_{n-1}(U_{\b_p^\f},p).
$$
\end{coro} 

These two results establish that in these cases, if 
$U_{\b_p^\f}$ is connected, then $\pi_1(\of,\f_{x_0})$ is trivial. This is 
granted if $p\b p$ is a von Neumann algebra. However, note that  
$\pi_2(\of,\f_{x_0})$ is not trivial. 
This is because $\b_p^\f$ is a finite von Neumann algebra and 
therefore $U_{\b_p^\f}$ has non trivial fundamental group 
(see \cite{han}, \cite{sch}).

\section{Vector states in $\l$}
In \ref{soporte} it was shown that a state $\Phi$ of $\l$ with support 
$e=e_x$ for some $x \in \sp$ is of the form $\Phi=\f_x$ for some state $\f$ in 
$\b$ with support $p$. Recall that we denote by $\Sigma_p(\b)$ the set of 
states of $\b$ with support $p$, and by $\Sigma_{p,X}$ the set of states of 
$\l$ with support equivalent to $e$. In other words, 
$\displaystyle{\Sigma_{p,X}=\cup_{\f \in \Sigma_p(\b)}\of}$. One has the assignment
$$
\sp \times \Sigma_p(\b) \to  \Sigma_{p,X} \quad ,\quad (x,\f) \mapsto \f_x.
$$

Remember that $\f_x=\psi_y$, with $\f,\psi \in \Sigma_p(\b)$, $x,y \in \sp$ 
if and only if $\psi=\f \circ Ad(u)$ with $y=xu$ and $u\in U_{p\b p}$ (see
\ref{fibra de fix} part c)).
 
The unitary group $U_{p\b p}$ acts both on $\sp$ and $\Sigma_p (\b)$. 
We may consider the diagonal action on $\sp \times \Sigma_p(\b)$, defined
by $u.(x,\f)=(xu,\f \circ Ad(u))$. It follows that if we denote the quotient
$$
\sp \times \Sigma_p(\b)/ \{(x,\f)\sim (xu,\f \circ Ad(u)) , u\in U_{p\b p}\} 
:= \sp \times_{U_{p\b p}} \Sigma_p(\b)
$$
(as is usual notation), then the assignment above induces a bijection
$$
\sp \times_{U_{p\b p}} \Sigma_p(\b)\simeq \Sigma_{p,X}.
$$
If we endow $\sp \times_{U_{p\b p}} \Sigma_p(\b)$ with the quotient topology 
(where $\sp$ and $\Sigma_p(\b)$ are considered with the norm topologies), 
a natural question is what topology does this bijection induce in 
$\Sigma_{p,X}$. The following result states that convergence of a sequence in 
the quotient topology is equivalent in $\Sigma_{p,X}$ to convergence (in norm) 
of the states and their supports.
\begin{prop}
Consider in $\Sigma_{p,X}$ the metric $d$ given
by 
$$
d(\Phi, \Psi)=\|\Phi - \Psi \| + \|\supp(\Phi)-\supp(\Psi) \|.
$$
Then the metric space $(\Sigma_{p,X},d)$ is homeomorphic to 
$\sp \times_{U_{p\b p}} \Sigma_p(\b)$, where the homeomorphism is given by the 
above bijection.
\end{prop}
\begin{proof}{Proof.}
Denote by $[(x,\f)]$ the class of $(x,\f)$ in $\sp \times_{U_{p\b p}} \Sigma_p(\b)$. 
Suppose that $[(x_n,\f_n)]$ converge to $[(x,\f)]$ in 
$\sp \times_{U_{p\b p}} \Sigma_p(\b)$. Then there exist unitaries 
$u_n$ in $p\b p$ such that $x_nu_n$ tends to $x$ and $\f_n\circ Ad(u_n)$ tends 
to $\f$ in the respective norms. By continuity of the inner product, it is clear 
then that $e_{x_n} =\theta_{x_n,x_n}=\theta_{x_n u_n,x_n u_n} \to e_x$ and 
${\f_n}_{x_n}=(\f_n\circ Ad(u_n))_{x_n u_n}\to \f_x$, and therefore the 
assignment $[(x,\f)]\mapsto \f_x$ is continuous. On the other direction,
suppose that $d(\Phi_n,\Phi)$ tends to zero. There exist 
$\f_n, \f \in \Sigma_p(\b)$ and $x_n, x \in \sp$ such that $\Phi_n={\f_n}_{x_n}$ 
and $\Phi=\f_x$. We have that $\supp(\Phi_n)=e_{x_n} \to \supp(\Phi)=e_x$. Now 
$e_{x_n}=\rho(x_n)$, $e_x=\rho(x)$, and $\rho:\sp \to \e_e$ is a fibre bundle 
with fibre $U_{p\b p}$, therefore there exist unitaries $u_n$ in $p\b p$ such 
that $x_nu_n \to x$. We may replace the $x_n$ by $y_n=x_nu_n$ and $\f_n$ by 
$\psi_n=\f_n\circ Ad(u_n)$, and still have $\Phi_n={\psi_n}_{y_n}$, with 
$y_n\to x$. We claim that $\psi_n \to \f$. Indeed, if $a\in \b$, by a typical 
argument
$$
|\psi_n (a)-\f(a)|=|\Phi_n(\theta_{y_na,y_n}) -\Phi (\theta_{xa,x})|\le 
\|\Phi_n\|\|\theta_{y_na,y_n}-\theta_{xa,x}\|+
\|\Phi_n-\Phi\|\|\theta_{xa,x}\|.
$$
The first summand is bounded by (using the Cauchy-Schwarz inequality)
$$
\|\theta_{y_na,y_n}-\theta_{xa,x}\|\le \|\theta_{y_na,y_n-x}\|+
\|\theta_{y_na-xa,x}\|\le \|y_n\|\|a\|\|y_n-x\|+\|y_n-x\|\|a\|\|x\|,
$$
which equals $2\|a\|\|y_n-x\|$.
The other summand equals $\|\Phi_n-\Phi\|\|a\|$. It follows that 
$[(y_n,\psi_n)]=[(x_n,\f_n)]\to [(x,\f)]$.
\end{proof}
Since $\|\Phi -\Psi\|\le d(\Phi,\Psi)$, it is clear that the 
inclusion $(\Sigma_{p,X},d)\subset (\l^*, \| \ \|)$ is continuous.
The following example shows that the topology given by the 
metric $d$ in $\Sigma_{p,X}$ does not coincide with the norm
 topology of the conjugate space of $\l$. In other words, that convergence
  of the vector states (which a priori have equivalent supports) 
  does not imply convergence of the supports.
\begin{exam}\label{ejemplo}\rm
Let $\b=D\subset B(\ell^2(\zN))$ be the subalgebra of diagonal matrices (with respect 
to the canonical basis). Consider the conditional expectation $E: B(\ell^2(\zN))\to D$ 
which consists on taking the diagonal entries. Let $a\in D$ be a trace class positive diagonal
operator with trace one, and no zero entries in the diagonal. Put $\f(x)=Tr(ax)$, 
$x\in B(\ell^2(\zN))$. Clearly, $\f$ is faithful and $B(\ell^2(\zN))^\f=D$. Let $b$ 
be the unilateral shift in $\ell^2(\zN)$. Denote by $q_n$ the $n \times n$ Jordan 
nilpotent, and $w_n$ the unitary operator on $\ell^2(\zN)$ having the unitary 
matrix $q_n+q_n^{* \ n-1}$ on the first $n\times n$ corner and the rest of
the diagonal completed with $1$. 

We shall consider the $D$-right module $X$ as the completion of $B(\ell^2(\zN))$ 
with the $D$-valued inner product given by $E$, i.e. $\langle x,y\rangle =E(x^*y)$, 
$x,y \in  B(\ell^2(\zN))$. Note that $X$ is also a $B(\ell^2(\zN))$-left module, 
and so the elements of $B(\ell^2(\zN))$ act as adjointable operators in $X$.  
Consider the faithful state $\f$ of $D$ equal to the restriction of the former $\f$. 
The elements $w_n$ and $b$ lie in the unit sphere of $X$. 

We claim that the projections $e_{w_n}$ do not converge to $e_b$. Suppose that they 
do converge, again using that the map $x\mapsto e_x$ is a bundle, there would exist 
unitaries $v_n\in D$ such that $w_nv_n \to b$ in $\s1$. 
In \cite{av2} it was shown that the element $b$ cannot be approximated by unitaries 
of $B(\ell^2(\zN))$ in the norm topology of the module $X$.

On the other hand, the states $\f_{w_n}$ converge to $\f_b$ in the 
norm topology of the conjugate space of ${\cal L}_\b(X)$. Indeed
$$
|\f_{w_n}(t)-\f_b(t)|=\left|Tr\left(a(\langle w_n,t(w_n)\rangle -\langle b,t(b)\rangle )\right)\right| 
$$
$$
\le \left|Tr\left(a(\langle w_n,t(w_n)-t(b)\rangle )\right)\right|+
\left|Tr\left(a(\langle w_n-b,t(b)\rangle )\right)\right|.
$$
The first summand can be bounded by $\|t\|\ Tr(a(2-E(w_n^*b)-E(b^*w_n)))$. 
Since $Tr(a)=1$ and $E$ is trace invariant, this term equals
$$
\|t\|\ \left(2-Tr(aw_n^*b+ab^*w_n)\right)=2 \|t\|\   \sum_{k\ge n} a_k,
$$
where $a_k$ are the diagonal entries of $a$. It is clear that this term tends to zero
when $n\to \infty$. 
The other summand can be dealt in a similar way, establishing our claim.

Summarizing, the states $\f_{w_n}$ converge but their supports $e_{w_n}$ do not.

\end{exam}
Next we shall see that the quotient map
$$
\wp_1:\sp \times \Sigma_p(\b)\to  \sp \times_{U_{p\b p}} \Sigma_p(\b) \quad , 
\quad  \wp_1(x,\f)=[(x,\f)]
$$
and the projection
$$
\wp_2:\sp \times_{U_{p\b p}} \Sigma_p(\b) \to \sp / U_{p\b p} \quad , 
\quad \wp_2([(x,\f)])=[x]
$$
are fibrations. Equivalently, if $\Sigma_{p,X}$ is considered with the topology 
induced by the metric $d$, the maps $(x,\f) \mapsto \f_x$ and $\f_x \mapsto e_x$ 
are fibrations.  In what follows, for brevity, we shall use  $\Sigma_{p,X}$ 
(considered with the metric $d$) 
instead of $\sp \times_{U_{p\b p}} \Sigma_p(\b)$. Therefore $\wp_1(x,\f)=\f_x$.
\begin{theo}
The map $\wp_1 :\sp \times \Sigma_p(\b) \to \Sigma_{p,X}$, 
$\wp_1 (x,\f)=\f_x$ is a principal fibre bundle with fibre $U_{p\b p}$.
\end{theo}
\begin{proof}{Proof.}
It suffices to exhibit a local cross section around a generic base point $\f_x$. We 
claim that there is a neighborhood of $\f_x$ such that elements $\psi_y$ 
in this neighborhood satisfy that $\langle y,x\rangle $ is invertible. Indeed, if $d(\f_x,\psi_y)< r$, 
then $\|e_x -e_y\|< r$. If we choose $r$ small enough so that $e_y$ lies in the ball 
around $e_x$ in which a local cross section of $\rho (x)=e_x$ is defined, then there 
exists a unitary $u$ in $p\b p$ such that $\|x-yu\|< 1$. Note that 
$$
\|p-\langle yu,x\rangle \|=\|\langle x-yu,x\rangle \|\le \|x-yu\|< 1.
$$
Then $\langle yu,x\rangle =u^*\langle y,x\rangle $ is invertible in $p\b p$, and therefore also $\langle y,x\rangle $. 
In this neighborhood put 
$$
s(\psi_y)=(y \mu (\langle y,x\rangle ),\psi \circ Ad(\mu(\langle y,x\rangle )),
$$ 
where $\mu$ denotes the unitary part in the polar decomposition of invertible elements in
 $p\b p$ as before. We claim that $s$ is well defined, is a local cross section and is 
 continuous. 

Suppose that $\psi_y=\psi'_{y'}$, then there exits a unitary $u$ in $p\b p$ such 
that $y'=yu$ and $\psi=\psi'\circ Ad(u)$. Then $y'\mu(\langle y',x\rangle )=yu \mu(u^*\langle y,x\rangle )=y\mu(\langle y,x\rangle )$. 
Also, 
$\psi'\circ Ad(\mu(\langle y',x\rangle ))=\psi \circ Ad(u) \circ Ad(\mu(u^*\langle y,x\rangle ))=
\psi \circ Ad(u)\circ Ad(u^*) \circ Ad(\mu(\langle y,x\rangle ))=\psi \circ Ad(\mu(\langle y,x\rangle ))$.

That it is a cross section is apparent. Let us see that $s$ is continuous. Suppose that 
${\psi_n}_{y_n} \to \f'_{x'}$ for $\f'_{x'}$ in the neighborhood of $\f_x$ 
where $s$ is defined. This implies that there exist unitaries $u_n$ in $U_{p\b p}$ 
such that $y_nu_n \to x'$ and $\psi_n \circ Ad(u_n) \to \f'$ in the norm topologies. 
The continuity of the inner product implies that 
$y_nu_n \mu(\langle y_nu_n,x\rangle )=y_n\mu(\langle y_n,x\rangle )\to x' \mu(\langle x',x\rangle )$. Also 
$\psi_n \circ Ad(u_n) \circ Ad(\mu(\langle y_nu_n,x\rangle ))=
\psi_n \circ Ad(\mu(\langle y_n,x\rangle )) \to \f' \circ Ad(\mu(\langle x',x\rangle ))$.
\end{proof}
Next we consider the map $\wp_2 :\Sigma_{p,X}\to \sp / U_{p\b p}$. 
Recall that $\sp / U_{p\b p}\simeq \e_e$, where the homeomorphism is given by 
$[x]\mapsto e_x$. The following result states that taking support of a state in 
$\Sigma_{p,X}$ (regarded with the $d$ topology) is a fibration.

\begin{theo}\label{wp2}
The map $\wp_2 :\Sigma_{p,X} \to \sp / U_{p\b p}$, given by
$\wp_2(\f_x)=[x]$ is a fibration with fibre $\Sigma_p(\b)$.
\end{theo}
\begin{proof}{Proof.}
Consider the diagram 
$$
\diagram{
\sp \times \Sigma_p(\b) &  \hrarr^{\wp_1}_{} & \Sigma_{p,X}  \cr
{} & \ddrarr{p }{}       & \vdarr{}{\wp_2}\cr
{}  &      {}        &  \sp / U_{p\b p} ,
}
$$
where $p$ is given by $p(x,\f)=[x]$. Clearly $p$ is a fibre bundle, because it is the composition
of the projective bundle $x \mapsto [x]$ with the projection $(x,\f)\mapsto x$. The map $\wp_1$ was shown to be a 
fibration.
It follows from \ref{lemita} the $\wp_2$ is a fibration. 
The fibre $\wp_2^{-1}([x])$ consists of all states $\f_y$ with $[y]=[x]$. 
Then there exists $u\in U_{p\b p}$ such that $\f_y=(\f\circ Ad(u^*))_x$, so that one
may fix $x$ (and not just $[x]$). Now $\f_x=\psi_x$ implies $\f=\psi$. It follows that the fibre 
over $[x]$ is the set $\{\f_x: \f \in \Sigma_p(\b)\}$, which identifies with $\Sigma_p(\b)$.
\end{proof}

We will use the fibrations $\wp_1$ and $\wp_2$ to obtain information about the homotopy 
type of these spaces.

As in the previous section, applying the homotopy exact sequences 
of these fibrations, one obtains
%\fontsize {12}{16}\selectfont
$$
\dots \pi_n\left(U_{p\b p},p\right) \to \pi_n(\sp \times \Sigma_p(\b),(x_0,\f)) 
\stackrel{(\wp_1)_*}{\to}  
 \pi_n(\Sigma_{p,X},\f_{x_0}) \to \pi_{n-1}(U_{p\b p},p)\to \dots
$$
%\fontsize {14}{16}\selectfont
and
$$
\dots \pi_n(\Sigma_p(\b),\f)\to \pi_n(\Sigma_{p,X},\f_{x_0}) 
\stackrel{(\wp_2)_*}{\to}
\pi_n({\cal E},e_{x_0})\to \pi_{n-1}(\Sigma_p(\b),\f)\dots
$$

First note that since $\Sigma_p(\b)$ is convex, $\sp \times \Sigma_p(\b)$ 
has the same homotopy type as $\sp$, and 
$$
\pi_*(\Sigma_{p,X})=\pi_*({\cal E}_{e_x}).
$$

\begin{coro}
The space ${\cal S}_p(H_\b) \times  \Sigma_p(\b)$  is contractible.
\end{coro}
\begin{proof}{Proof.}
It was remarked in the preceding section that ${\cal S}_p(H_\b)$ is contractible.
\end{proof}
\begin{coro}
For  $\f_0 \in \Sigma_p(\b)$ and $x_0 \in {\cal S}_p(H_\b)$ fixed,
$$
\pi_n\left(\Sigma_{p,H_\b},{\f_0}_{x_0}\right)\cong 
\pi_{n-1}(U_{p\b p},p), \ \ n\ge 1.
$$
In particular, if $U_{p\b p}$ is connected, 
$\pi_1\left(\Sigma_{p,H_\b},{\f_0}_{x_0}\right)=0$. 
If moreover $p\b p$ is a properly infinite von Neumann algebra, 
$\Sigma_{p,H_\b}$ has trivial homotopy groups of all 
orders.
\end{coro}
\begin{proof}{Proof.}
The first fact follows from the contractibility of ${\cal S}_p(H_\b) \times  \Sigma_p(\b)$, 
which implies that in the homotopy sequence 
$\pi_k({\cal S}_p(H_\b) \times\Sigma_p(\b),(\f_0,x_0))=0$ for all $k$. 
The second fact follows using that $U_{p\b p}$ is connected. Using that (\cite{bru}) if 
$p\b p$ is properly infinite, then $U_{p\b p}$ is contractible, it follows that 
$$
\pi_n(\Sigma_{p,H_\b},{\f_0}_{x_0})=0
$$
for all $n\ge 0$. 
\end{proof}
One can be more specific, since the homotopy groups of the unitary group of a 
C$^*$-algebra (at least for $n=1$) have been computed in many cases (\cite{han}, 
\cite{sch}, \cite{zha}). For example, in the von Neumann algebra case, one can 
compute the fundamental group of the unitary group in terms of the type decomposition 
of the algebra.
\begin{coro}
If $p\b p$ is a finite von Neumann algebra, then $\Sigma_{p,X}$ is connected.
\end{coro}
\begin{proof}{Proof.}
It was noted before that if $p\b p$ is finite, then $\sp$ is connected.
\end{proof}
\begin{coro}
If $p\b p$ is a properly infinite 
%von Neumann 
algebra,  then for $n\ge 0$
$$
\pi_n\left(\Sigma_{p,X},{\f_0}_{x_0}\right)
\cong \pi_n(\sp, x_0).
$$
If moreover $Xp$ is selfdual, 
then 
$$
\pi_n\left(\Sigma_{p,X},{\f_0}_{x_0}\right)=0
$$
for all $n\ge 0$.
\end{coro}
\begin{proof}{Proof.}
The proof follows writing the homotopy exact sequence of $\wp_1$. If $p\b p$ is 
properly infinite, its unitary group is contractible. If moreover $Xp$ is selfdual, 
it was pointed out before that $\sp$ is contractible.
\end{proof}
We turn now our attention to the bundle $\wp_2$.

\begin{coro} \label{teorema}
If $p\b p$ is a von Neumann algebra and $Xp$ is selfdual, then 
$\pi_1(\Sigma_{p,X},\f_x)=0$

\end{coro}

\begin{proof}{Proof.}
It was shown in \cite{acs2} that $\pi_1({\cal E}_{e_x})=0$
\end{proof}

\begin{rema} \rm
There is another map related to this situation, namely the other projection $\wp_3$,
$$
\wp_3 :\Sigma_{p,X} \to \Sigma_p(\b)/U_{p\b p} ,\ \wp_3(\f_x)=[\f].
$$
This map is well defined and continuous, if one goes back to the notation 
$\sp \times_{u_{p\b p}} \Sigma_p(\b)$, $\wp_3$ is the map $(x,\f)\mapsto \f$ 
at the quotient level,
$$
[(x,\f)]\mapsto [\f].
$$
However this map is not, in general, even a weak fibration. To see this consider the case when $X=\b$ 
is a finite algebra, and $p=1$. Here ${\cal L}_\b(\b)=\b$ and $\Sigma_{1,\b}$ consist of the 
states of $\b$ with support equivalent to 1 (note that $x\in \s1$ verifies $x^*x=1$, i.e. $x \in U_\b$, 
and $e_x=1$). That is, $\Sigma_{1,\b}$ is the set of faithful states of $\b$ 
($=\Sigma_1(\b)$ in our notation). It follows that $\wp_3$ is just the quotient map
$$
\Sigma_1(\b) \to \Sigma_1(\b)/U_\b.
$$
Moreover, take $\b=M_n(\zC)$ ($n< \infty$). Then the quotient map above
is not a weak fibration. 
Indeed,
both sets $\Sigma_1(M_n(\zC))$ and $\Sigma_1(M_n(\zC))/U_{M_n(\zC)}$ are convex metric spaces. 
The latter can be identified, using the density matrices, as the $n$-tuples of eigenvalues
$(\lambda_1,...,\lambda_n)$ arranged in decreasing order and normalized such that 
$\sum \lambda_k=1$, with the $\ell_1$ distance. If this quotient map were a weak fibration, 
then the fibre would have trivial homotopy groups of all orders. 
This is clearly not the case, since the fibre is the unitary group $U(n)$ of $M_n(\zC)$.
\end{rema}

\begin{rema} \rm
The set $\of$ lies inside $\Sigma_{p,X}$, namely as the states $\f_x$ with $\f$ fixed. 
If one regards $\of$ with the metric $d_\f$ and $\Sigma_{p,X}$ with the metric $d$, 
it is clear that the inclusion is continuous. Indeed,
it was noted that $\supp$ is continuous in $\of$. Therefore if $d_\f(\f_{x_n},\f_x)\to 0$, 
then $e_{x_n}\to e_x$, which implies that 
$d(\f_{x_n},\f_x)\to 0$.

However, the identity mapping $(\of,d_\f)\to (\of,d)$ is not (in general) a 
homeomorphism. 
Indeed, take  $X=\b$ and  $\f$ faithful. Then $\of$ is the unitary orbit $\{\f\circ Ad(u): u\in U_\b\}\sim U_\b/U_{\b^\f}$, and it is clear that $d_\f$ induces the same topology as the quotient topology ($U_\b$ with the norm topology). On the other hand, $\Sigma_{p,X}$ coincides in this case with $\Sigma_1(\b)$ the set of faithful states of $\b$, and the metric $d$ is just the usual norm of the conjugate space $\b^*$. In \cite{av1} it was shown that in general, the unitary orbit does not have norm continuous local cross sections to the unitary group, but it does have local cross sections which are continuous in the quotient topology $U_\b/U_{\b^\f}$.
\end{rema}

\begin{rema} \rm
The metric $d(\Phi,\Psi)=\|\Phi-\Psi\|+\|\supp(\Phi)-\supp(\Psi)\|$ in the state space of 
$\l$ is weird. For example, in this metric, $\Sigma_{p,X}$ is open. Moreover, any state $\psi$ of $\l$ such that $d(\Psi,\Sigma_{p,X})< 1$, actually lies in $\Sigma_{p,X}$. 
Indeed, if $\Phi$ is 
a state of $\l$, and $d(\Phi, \f_x)< 1$ for some $x\in \sp$ and $\f \in \Sigma_p(\b)$, 
then $\|supp(\Phi)-e_x\|< 1$, and therefore $\supp(\Phi)$ and $e_x$ are unitarily 
equivalent. That is $\supp(\Phi)=e_y$, with $y=U(x)$ for some unitary $U$ in $\l$. 
Then, by \ref{soporte}, there exists $\psi \in \Sigma_p(\b)$ such that 
$\Phi=\psi_y \in \Sigma_{p,X}$. 

However, if $X=\b$ and $\b$ is finite dimensional, then the topology of the $d$-metric coincides in $\Sigma_{p,\b}$ with 
the usual norm topology. Indeed, it suffices to see that the map $\f_x \mapsto e_x$ is continuous in the
norm topology. Since we are in the finite dimensional case, it suffices to argue with (positive) density
matrices, with trace 1. Note that the states of the form $\f_x$ have equivalent supports, i.e. their
density matrices have kernels with the same dimension. Suppose that $a_n$ is a sequence of positive matrices
with trace 1 and $nul(a_n)=k$, converging in norm to the matrix $a$, also with (a priori) $nul(a)=k$. Then the
projections $P_{\ker a_n}$ onto the kernels converge in norm to $P_{\ker a}$. Indeed, we claim that 
one can find an open interval around zero and an integer $n_0$ such that for $n\ge n_0$ no 
eigenvalue of $a_n$ (other than zero) lies inside this interval. And by a routine spectral theory (Riesz integral) argument, one has that
 $P_{\ker a_n} \to P_{\ker a}$. If one could find no such interval, then there
would exist a sequence $\lambda_n > 0$ such that $\lambda_n$ is an eigenvalue of $a_n$ and $\lambda_n \to 0$.
If $q_n$ is the spectral projection corresponding to $\lambda_n$, then $a_n=b_n+\lambda_n q_n$. Then 
$b_n \to a$, where $nul (a)=k$ and $nul (b_n)< k$, which cannot happen.
\end{rema}

\section{Purification of $\Sigma_{p,X}$}

In this section $\b$ is a von Neumann algebra. We will consider a natural representation for $\l$, studied in \cite{pas} and \cite{rie}, in which all the states $\f_x \in \Sigma_{p,X}$ are induced by vectors in this Hilbert space.
Consider the algebraic tensor product $X\otimes H$, where $H$ is a Hilbert space on which $\b$ acts. We will choose $H$ the space of a standard representation of $\b$. Recall the fact that for such a representation there exists a cone ${\cal P}$, called the positive standard cone, with many remarkable properties. Among them, any positive normal functional in $\b$ is implemented by a unique vector in this cone. In the vector space $X\otimes H$ consider the semidefinite positive form given by $[x\otimes \xi, y\otimes \eta]=(\xi, \langle x,y\rangle \eta)$, where $(\ , \ )$ is the inner product of $H$. Denote by $Z=\{z\in X\otimes H: [z,z]=0\}$, and let ${\cal H}$ be the Hilbert space obtained as the completion of $X\otimes H/Z$. The representation $\rho:\l \to B({\cal H})$ is given by $\rho (t)([x\otimes \xi])=[t(x)\otimes \xi]$. 

\begin{lemm}
In the representation $\rho$, the state $\f_x \in \Sigma_{p,X}$ is implemented by the vector $x\otimes \xi$, where $\xi$ is the unique vector in the positive cone ${\cal P}$ which implements $\f$ ($\f(a)=\omega_\xi(a)=(a\xi,\xi)$), that is
$$\f_x(t)=[\rho(t)(x\otimes \xi), x\otimes \xi], \ t\in \l.$$
\end{lemm}
\begin{proof} {Proof.}
Straightforward: $[\rho(t)(x\otimes \xi),x\otimes \xi]=(\xi,\langle t(x),x\rangle \xi)=(\langle x,t(x)\rangle \xi,\xi)=\f(\langle x,t(x)\rangle )=\f_x(t)$.
\end{proof}

In order to simplify the exposition, we shall restrict to the case $p=1$. This is in fact no restriction at all, since the general case can be easily reduced to this situation 
(note that $\sp$ is the unit sphere of the $p\b p$ module $Xp$). 
Therefore, the unitary vectors of the cone implementing the faithful states of $\b$ are the vectors which are cyclic and separating for $\b$. Let us denote ${\cal A}(X)=\{[x\otimes \xi]: x\in {\cal S}_1(X), \xi \in {\cal P}, \hbox{ cyclic and separating for }\b, \|\xi\|=1\}\subset {\cal H}$.
\begin{lemm}
Let $x,y \in \s1$ and $\xi,\eta \in {\cal P}$ unit, cyclic and separating, then the elements $x\otimes \xi$ and $y\otimes \eta$  induce the same element in ${\cal A}(X)$ only if $x=y$ and $\xi=\eta$. In other words, there is a bijection
$$
\s1 \times \Sigma_1(\b) \leftrightarrow {\cal A}(X), \ (x,\f)\mapsto x\otimes \xi .
$$
\end{lemm}
\begin{proof}{Proof.}
Suppose that $x\otimes\xi \sim y\otimes \eta$, with $x,y,\xi,\eta$ as above. Then 
$$
0=[x\otimes\xi - y\otimes \eta,x\otimes\xi - y\otimes \eta]=2-2Re((\xi,\langle x,y\rangle \eta)).$$
That is, $(\xi,\langle x,y\rangle \eta)=1$. Since $\xi, \eta$ are unital and $\|\langle x,y\rangle \|\le 1$, by the Cauchy-Schwarz inequality this implies that $\langle x,y\rangle \eta=\lambda \xi$, for $\lambda \in \zC$, $|\lambda|=1$. Again, using that $\xi, \eta$ are unital, this implies $\lambda=1$, i.e. $\langle x,y\rangle \eta=\xi$.
On the other hand, the states induced in $\l$ by the vectors $[x\otimes\xi]$ and $[y\otimes \eta]$ via the representation $\rho$ were shown to be $\f_x$ and $\psi_y$, where $\f, \psi$ are the states induced in $\b$ by $\xi,\eta$, respectively, as shown in the lemma above. By \ref{fibra de fix}, $\f_x=\psi_y$ implies that there exists $u\in U_\b$ such that $y=xu$ and $\psi=\f\circ Ad(u)$. So $\langle x,y\rangle \eta=\xi$ translates into $u\eta=\xi$. The other identity $\psi=\f\circ Ad(u)$ can also be interpreted in terms of these vectors in the cone ${\cal P}$. Namely, the unique vector in the cone associated to the state $\f\circ Ad(u)$ is $u^*Ju^*J\xi$, where $J$ denotes the modular conjugation of the standard representation $\b \subset B(H)$. Indeed, clearly $u^*Ju^*J\xi \in {\cal P}$, and
$(au^*Ju^*J\xi, u^*Ju^*J\xi)=(Ju^*Jau^*\xi, Ju^*Ju^*\xi)=(au^*\xi,u^*\xi)=(uau^*\xi,\xi)=\f(uau^*)$. Therefore, by the uniqueness condition (on vectors in the cone inducing states), it follows that $\eta=u^*Ju^*J\xi$. Combining this with $u\eta=\xi$ yields
$$
\xi=Ju^*J\xi=Ju^*\xi, \ \hbox{ i.e. } \xi=u^*\xi.$$
This implies that $u^*$ acts as the identity operator on $\b '\xi$, which is dense in $H$, because $\xi$ is cyclic for $\b '$. Therefore $u=1$. Then $x=y$ and $\xi=\eta$. 
\end{proof}

These two lemmas state that
$$
\wp_1: \sp\times \Sigma_1(\b)\to \Sigma_{1,X} , \ \wp_1(x,\f)=\f_x ,
$$
in this representation looks like
$$
\vec{\wp_1}: {\cal A}(X)\to \Omega_{{\cal A}(X)} , \ \vec{\wp_1}([x\otimes\xi])=\omega_{[x\otimes\xi]},
$$
where $\omega_{[x\otimes\xi]}$ is the state induced by the vector $[x\otimes\xi] \in {\cal H}$, and $\Omega_{{\cal A}(X)}$ is the space of all such states with symbols in ${\cal A}(X)$. What one gains by taking this standpoint is that ${\cal A}(X)$ has a natural topology, as a subset of the Hilbert space ${\cal H}$. The set $\Omega_{{\cal A}(X)}\sim \Sigma_{1,X}$ is therefore endowed with a weaker topology than the $d$ metric, namely the quotient topology induced by ${\cal A}(X)$ and $\vec{\wp_1}$. The fibre of this map is , as before, a copy of the unitary group $U_\b$ of $\b$. The next result examines how the unitary group $U_\b$ appears inside ${\cal A}(X)$ and which is its relative topology. By the above result, we can  omit the brackets when dealing with classes of elementary tensors of the form $x\otimes\xi$ in ${\cal A}(X)\subset {\cal H}$ ($x\in\s1$, $\xi$ unit, cyclic and separating in ${\cal P}$). Also, note that the vectors $x\otimes \xi \in {\cal A}(X)$ are cyclic for $\rho$, but not separating in general.

\begin{prop}
Given a fixed element $x\otimes \xi \in {\cal A}(X)$, the fibre $\vec{\wp_1}^{-1}(\omega_{x\otimes\xi})$
is the set $\{xu\otimes u^*Ju^*J\xi: u\in U_\b\}$ which is in one to one correspondence with $U_\b$. The relative topology induced on $U_\b$ by this bijection is the strong operator topology.
\end{prop}
\begin{proof}{Proof.}
If $y\otimes\eta$ lies in the fibre $\vec{\wp_1}^{-1}(\omega_{x\otimes\xi})$, then $\omega_{x\otimes\xi}=\omega_{y\otimes\eta}$, or $\f_x=\psi_y$, where as in the previous lemma $\f$ and $\psi$ are the states of $\b$ associated to the vectors $\xi$ and $\eta$. Again, this implies that there exists a unitary in $U_\b$ such that $y=xu$ and $\eta=u^*Ju^*J\xi$. Then
$y\otimes \eta=xu\otimes u^*Ju^*J\xi$.
 Now suppose that a net $xu_\alpha\otimes u_\alpha^*Ju_\alpha^*J\xi$ converges to $xu\otimes u^*Ju^*J\xi$ in the Hilbert space topology (of ${\cal H}$). This implies that 
 $$
 \|xu_\alpha\otimes u_\alpha^*Ju_\alpha^*J\xi - xu\otimes u^*Ju^*J\xi\|^2=2-2Re((u_\alpha^*Ju_\alpha^*J\xi,\langle xu_\alpha,xu\rangle u^*Ju^*J\xi))\to 0 
 %\hbox{ with } \alpha .
 $$
 with $ \alpha$.  In other words, 
 $$
 (u_\alpha^*Ju_\alpha^*J\xi,u_\alpha^*Ju^*J\xi)=(u^*\xi,u^*\xi)\to 1 .
 $$
 Equivalently, $(uu_\alpha^*\xi,\xi)\to 1$. This implies that $\|(u-u_\alpha)\xi\|\to 0$ in the Hilbert space norm (of $H$). Now let $a'\in \b '$, then 
 $$
 \|(u-u_\alpha)a'\xi\|=\|a'(u-u_\alpha)\xi\|\le \|a'\|\|(u-u_\alpha)\xi\|\to 0 .$$
 That is $u_\alpha \nu \to u\nu$ in a dense subset of vectors $\nu \in H$. Since $u,u_\alpha$, being unitaries, are bounded in norm, this implies strong operator convergence of $u_\alpha$ to $u$. The converse implication is straightforward.
 \end{proof}
The tensor product $X\otimes H/Z$ is a $\b$-bimodule tensor product, in the sense that for any $b\in \b$, $x\in X$ and $\nu \in H$, one has $xb\otimes \nu$ equivalent to $x\otimes b\nu$.
Then the elements $xu\otimes u^*Ju^*J\xi$ in the fibre of $\omega_{x\otimes \xi}=\f_x$ can be parametrized $x\otimes Ju^*J\xi=x\otimes Ju^*\xi$. We prefer the first presentation because the vector $Ju^*\xi$ does not belong to ${\cal P}$. However the latter clarifies the action of $U_\b$ on ${\cal A}(X)$. Namely, the right-action 
$$
 (x\otimes \xi)\bullet u=x\otimes Ju^*J\xi .
$$
Note that it is indeed a right action:
$ (x\otimes \xi)\bullet vu=x\otimes J(vu)^*J\xi=x\otimes Ju^*JJv^*J\xi=((x\otimes \xi)\bullet v)\bullet u$. 

The sphere $\s1$ and the set $\Sigma_1(\b)$ of faithful states of $\b$ lie inside ${\cal A}(X)$. Pick a fixed element $x_0\in \s1$ and  $\xi_0 \in {\cal P}$ unit, cyclic and separating, inducing the state $\f_0$. The following maps are one to one:
$$
\s1 \to \{x\otimes \xi_0: x\in \s1 \}\subset {\cal A}(X) ,\ \   x\mapsto x\otimes\xi_0 ,
$$
and
$$
\Sigma_1(\b) \to \{x_0\otimes \xi: \xi \in {\cal P} \hbox{ unit, cyclic and separating} \},
$$
$$
\f \mapsto x_0\otimes \xi,
$$
where $\xi$ is the vector in the cone associated to $\f$.
\begin{prop}
The first bijection endows $\s1$ with the relative topology induced from ${\cal H}$, which is given by the following:
a net $x_\alpha$ converges to $x$ if and only if $\f_0(\langle x_\alpha -x,x_\alpha -x\rangle )\to 0$, if and only if $|x_\alpha -x|\to 0$ in the strong operator topology of $\b\subset B(H)$. The sphere $\s1 \subset X$ is closed in this topology.

The second bijection is a homeomorphism when $\Sigma_1(\b)$ is regarded with the norm topology and $\{x_0\otimes \xi: \xi \in {\cal P} \hbox{ unit, cyclic and separating} \}\subset {\cal H}$ is regarded with the Hilbert space norm of ${\cal H}$.
\end{prop}
\begin{proof}{Proof.}
The second statement is straightforward, because $\|x_0\otimes \xi-x_0\otimes\eta\|^2=2-2Re(\xi,\eta)=\|\xi-\eta\|^2$ and the well known fact that the topology of the distance between the vectors in ${\cal P}$ yields a topology which is equivalent to the one given by the norm of the induced states in the conjugate space.
Let  $x_\alpha\otimes \xi_0$ be a net, and $x\otimes \xi_0$ an element in ${\cal A}(X)$. Then $\|x_\alpha\otimes \xi_0-x\otimes \xi_0\|^2=2-2Re(\xi_0, \langle x_\alpha,x\rangle \xi_0) = 2-2Re(\f_0(\langle x_\alpha,x\rangle ))=\f_0(\langle x_\alpha-x,x_\alpha-x\rangle )$. Next we check that the convergence of the net in the sense described is equivalent to convergence to zero of $|x_\alpha -x|$ in the strong topology, where as is usual notation, $|y|=\langle y,y\rangle ^{1/2}$, for $y\in X$. Since $\f_0$ is implemented by the vector $\xi_0$, $\f_0(\langle x_\alpha -x,x_\alpha -x\rangle )=\| |x_\alpha -x|\xi_0\|^2$, and therefore convergence in the strong topology implies convergence in the former sense.
Suppose now that $\| |x_\alpha -x|\xi_0\|\to 0$, and take $a'\in \b '$. Then $\| |x_\alpha -x|a'\xi_0\|=\|a'|x_\alpha -x|\xi_0\|\le \|a'\|\| |x_\alpha -x|\xi_0\|\to 0$. The set $\{a'\xi_0: a'\in \b '\}$ is dense in $H$, and the operators $|x_\alpha -x|$ have bounded norms, therefore $|x_\alpha -x|$ tends strongly to zero.

Let us  prove now that the sphere $\s1\subset X$ is closed in this topology. First note that this topology, on norm bounded sets, is induced by the seminorms $n_\nu (x)=(\langle x,x\rangle \nu, \nu)$, $\nu \in H$, $\| \nu\|=1$ \cite{pas}.
Then it suffices to see that if $x_\alpha \to x$ with $x_\alpha \in \s1$, then $x\in \s1$.
Now, $(\langle x,x\rangle \nu,\nu)= 1$. Indeed, if $\omega_\nu(a)=(a\nu,\nu)$, then
$$
(\langle x_\alpha -x,x\rangle \nu,\nu)=\omega_\nu(\langle x_\alpha-x,x\rangle ) \le \omega_\nu(\langle x_\alpha-x,x_\alpha-x\rangle )^{1/2}=n_\nu(x_\alpha-x)^{1/2} ,
$$
i.e. $(\langle x_\alpha,x\rangle \nu,\nu)\to (\langle x,x\rangle \nu,\nu)$. Therefore
$$
0\leftarrow (\langle x_\alpha-x,x_\alpha-x\rangle \nu,\nu)=1+((\langle x,x\rangle -\langle x_\alpha,x\rangle -\langle x,x_\alpha\rangle )\nu,\nu)\to 1-(\langle x,x\rangle \nu,\nu).
$$
Since this is true for all unit vectors $\nu \in H$, it follows that $x\in \s1$.
\end{proof}
\begin{rema} \rm
Since $X$ is selfdual, it is a conjugate space \cite{pas}. The result above shows that the topology of $\s1$ induced by the Hilbert space norm of ${\cal H}$ coincides with the $w^*$ topology of $X\supset \s1$. Indeed, it was shown in \cite{pas} that a net $x_\alpha \to x$ in the $w^*$ topology if and only if $\varphi(\langle x_\alpha,y\rangle )\to \varphi(\langle x,y\rangle )$ for all $y\in X$, $\varphi \in \b_*^+$. This clearly implies that $\varphi(\langle x_\alpha-x,x_\alpha-x\rangle )\to 0$, which is the topology considered in the lemma (here the fact $\langle x,x\rangle =\langle x_\alpha,x_\alpha\rangle =1$ is crucial). Conversely
$$
\varphi(\langle x_\alpha-x,y\rangle )\le \varphi(\langle x_\alpha-x,x_\alpha-x\rangle )^{1/2}\varphi(\langle y,y\rangle )^{1/2}
$$
yields the other implication.

\end{rema}
We have examined the topologies induced on $\s1$ and $\Sigma_1(\b)$ by the described inclusions on ${\cal A}(X)$. We have seen before that ${\cal A}(X)\sim \s1 \times \Sigma_1(\b)$. These facts alone however do not imply that ${\cal A}(X)$ is homeomorphic to $\s1 \times \Sigma_1(\b)$ in the product topology (of the w$^*$ topology and the norm topology respectively). The next result shows that this is the case.
\begin{theo}
The bijection
$$
\s1 \times \Sigma_1(\b)\to {\cal A}(X)\ , \ (x,\f)\mapsto x\otimes \xi,
$$ 
is a homeomorphism when $\s1 \times \Sigma_1(\b)$ is endowed with the product topology of the w$^*$ topology of $\s1$ and the norm topology of $\Sigma_1(\b)$.
\end{theo}
\begin{proof}{Proof.}
By the result above, it is clear that if $x_\alpha \to x$ in $\s1$ and $\f_\beta \to \f$ in $\Sigma_1(\b)$, then $x_\alpha \otimes \xi_\beta \to x\otimes \xi$, where $\xi_\beta , \xi$ are the vectors in the positive cone inducing $\f_\beta, \f$. On the other direction, suppose that $x_\alpha \otimes \xi_\alpha \to x\otimes \xi$ in ${\cal A}(X)$. This means that $(\langle x,x_\alpha\rangle\xi_\alpha, \xi)\to 1$. Then, 
since $\|\langle x,x_\alpha\rangle\xi_\alpha\|\le 1$,
follows that
$\|\langle x,x_\alpha\rangle\xi_\alpha- \xi\|^2=
1+\|
\langle x,x_\alpha\rangle\xi_\alpha\|^2-
2Re(
\langle x,x_\alpha\rangle\xi_\alpha, \xi)
\to 0
$ and similarly 
$\|\langle x_\alpha,x\rangle\xi- \xi_\alpha\|\to 0$.
Then we get
\begin{equation}\label{vectores}
\langle x,x_\alpha\rangle\xi_\alpha- \xi\to 0
\ \mbox{ and }\
\langle x_\alpha,x\rangle\xi- \xi_\alpha\to 0
\end{equation}
This implies that
\begin{equation}\label{conv1}
\omega_{\langle x,x_\alpha\rangle\xi_\alpha}-\omega_\xi=
\f_\alpha\left(\langle x_\alpha,x\rangle\ \cdot\   \langle x,x_\alpha\rangle\right)-
\f(\ \cdot \ )
\longrightarrow 0
\end{equation}
Using \ref{vectores}  follows that  $J \langle x_\alpha,x\rangle J \langle x,x_\alpha\rangle\xi_\alpha-\xi_\alpha\to 0$
and so
\begin{equation}\label{conv2}
\omega_{J \langle x_\alpha,x\rangle J \langle x,x_\alpha\rangle\xi_\alpha}-\omega_{\xi_\alpha}=
\omega_{J \langle x_\alpha,x\rangle J \langle x,x_\alpha\rangle\xi_\alpha}-\f_\alpha\longrightarrow 0
\end{equation}
Note that for every $a\in \b$
\begin{eqnarray}\label{omegajota}
\omega_{J \langle x_\alpha,x\rangle J \langle x,x_\alpha\rangle\xi_\alpha}(a)& = &
(a J \langle x_\alpha,x\rangle J \langle x,x_\alpha\rangle\xi_\alpha,J \langle x_\alpha,x\rangle J \langle x,x_\alpha\rangle\xi_\alpha)\nonumber\\
 & = &
(\langle x_\alpha,x\rangle a\langle x,x_\alpha\rangle\xi_\alpha, J\langle x,x_\alpha\rangle\langle x_\alpha,x\rangle J\xi_\alpha).
\end{eqnarray}
But
\begin{eqnarray}
\|  J\langle x,x_\alpha\rangle\langle x_\alpha,x\rangle J\xi_\alpha -\xi_\alpha  \| &\leq &
\|J\langle x,x_\alpha\rangle\langle x_\alpha,x\rangle J\xi_\alpha-
J\langle x,x_\alpha\rangle\langle x_\alpha,x\rangle J\langle x_\alpha,x\rangle\xi\|
\nonumber\\
& &+ \|J\langle x,x_\alpha\rangle\langle x_\alpha,x\rangle J\langle x_\alpha,x\rangle\xi-
\xi_\alpha\|\nonumber
\end{eqnarray}
It is easy to prove that the first addend tends to zero using \ref{vectores}. The second one is equal to 
$ \|J\langle x,x_\alpha\rangle\langle x_\alpha,x\rangle J\langle x_\alpha,x\rangle J\xi-
\xi_\alpha\|=
\|JJ\langle x_\alpha,x\rangle J\langle x,x_\alpha\rangle\langle x_\alpha,x\rangle \xi-
\xi_\alpha\|$ which again tends to zero using \ref{vectores}.
Therefore, 
$  J\langle x,x_\alpha\rangle\langle x_\alpha,x\rangle J\xi_\alpha -\xi_\alpha \to 0$, and using \ref{omegajota} we get that
\begin{equation}\label{conv3}
\omega_{J \langle x_\alpha,x\rangle J \langle x,x_\alpha\rangle\xi_\alpha}-\f_\alpha
\left(\langle x_\alpha,x\rangle \cdot   \langle x,x_\alpha\rangle\right)
\longrightarrow 0
\end{equation}
Finally, using \ref{conv1}, \ref{conv2} and \ref{conv3} follows that $\f_\alpha\to \f$ in norm.

Equivalently, $\xi_\alpha \to \xi$ in $H$. Then $(\langle x,x_\alpha\rangle \xi_\alpha,\xi)\to 1$
implies that $(\langle x,x_\alpha\rangle \xi,\xi)\to 1$. Then $\f(|x-x_\alpha|)\to 0$,
i.e. $x_\alpha \to x$ in the w$^*$ topology.
\end{proof}
\begin{coro}
The space ${\cal A}(X)$ is homotopically equivalent to the sphere $\s1$ with the
w$^*$ topology.
\end{coro}
\begin{proof}{Proof.}
Recall that $\Sigma_1(\b)$ is convex.
\end{proof}
Now we focus on the map 
$$                                                              
\vec{\wp_1} : {\cal A}(X)\to \Omega_{{\cal A}(X)}, \vec{\wp_1}(x\otimes \xi)=\omega_{x\otimes \xi}.
$$
In order to see if this map is a fibration, we shall look for local cross sections. A powerful tool to state the existence of cross sections is Michael's theory of continuous selections \cite{mic}. An example of the use of this theory in the context of operator algebras is the paper by S. Popa and M. Takesaki \cite{pota}. To use Michael's theorem one must check first that the set function $\omega_{z\otimes \xi}\mapsto \vec{\wp_1}^{-1}(\omega_{z\otimes \xi})$ which assigns to each point in the base space the fibre over it, is {\it lower semicontinuous} \cite{mic}.

\begin{rema} \rm

In our context lower semicontinuity means that for any $r > 0$, and $x\otimes\xi$ the set $\{\omega_{y\otimes \eta}: \|y\otimes Ju^*\eta - x\otimes \xi\| < r\}$ is open in $\Omega_{{\cal A}(X)}$. In other words, for a state $\omega_{y\otimes\eta}$ close to $\omega_{x\otimes\xi}$ one should find an element $y\otimes Ju^*\eta$ in the fibre of $\omega_{y\otimes\eta}$ at distance less than $r$ to the fibre of $\omega_{x\otimes \xi}$. We have not specified yet the topology of this set $\Omega_{{\cal A}(X)}$. Lower semicontinuity implies that whatever topology one chooses, it must be stronger than the quotient topology  given by $\vec{\wp_1}$. Indeed, two states in $\Omega_{{\cal A}(X)}$ are close in this quotient topology if and only if there are elements in their fibres which are close in ${\cal A}(X)$. 

On the other hand this quotient topology is stronger than the norm topology. Recall Bures metric for states, defined as the infimum of the distances between vectors inducing the states taken over all possible representations where the two states are vector states. The topology induced by Bures metric on the state space is equivalent to the norm topology. This raises the question of whether this two topologies, the norm topology and the one induced by this purification coincide in $\Omega_{{\cal A}(X)} (=\Sigma_{1,X})$.

\end{rema} 

\begin{theo}
The quotient and the norm topology coincide in $\Omega_{{\cal A}(X)}.$
\end{theo}
\begin{proof}{Proof.}
It was noted that the quotient topology is stronger than the norm topology. Let us check the other implication. Let $\omega_{y_n \otimes \eta_n}$ be a sequence in $\Omega_{{\cal A}(X)}$ converging to $\omega_{x\otimes \xi}$ in norm. Testing convergence in operators of the form $\theta_{y_n a,y_n}$, $a\in \b$, yields
$$
\|a\| \|\omega_{y_n \otimes \eta_n}-\omega_{x\otimes\xi}\|=\|\theta_{y_na,y_n}\|\|\omega_{y_n \otimes \eta_n}-\omega_{x\otimes\xi}\|\ge |\omega_{y_n \otimes \eta_n}(\theta_{y_na,y_n})-\omega_{x\otimes\xi}(\theta_{y_na,y_n})|.
$$
Note that $\omega_{y_n \otimes \eta_n}(\theta_{y_na,y_n})=(a\eta_n,\eta_n)$ and $\omega_{x\otimes\xi}(\theta_{y_na,y_n})=(a\langle y_n,x\rangle \xi,\langle y_n,x\rangle \xi)$. This implies that 
$$
\|\omega_{\eta_n}-\omega_{\langle y_n,x\rangle \xi}\|\to 0 , \ \hbox{ as } n\to \infty.
$$
In particular, testing this difference at $1\in \b$, implies $(\langle x,y_n\rangle \langle y_n,x\rangle \xi,\xi)\to 1$. Therefore, 
$$
\|\langle x,y_n\rangle \langle y_n,x\rangle \xi-\xi\|^2=1+\|\langle x,y_n\rangle \langle y_n,x\rangle \xi\|^2-2Re(\langle x,y_n\rangle \langle y_n,x\rangle \xi,\xi)\to 0.
$$
Coming back to $\omega_{\eta_n}$ and $\omega_{\langle y_n,x\rangle \xi}$, note that the vectors $\eta_n$ belong to the cone ${\cal P}$, but not necesarilly the vectors $\langle y_n,x\rangle \xi$. However  $\delta_n=\langle y_n,x\rangle J\langle y_n,x\rangle \xi=\langle y_n,x\rangle J\langle y_n,x\rangle J\xi \in {\cal P}$ and we shall see that $\omega_{\delta_n}-\omega_{\langle y_n,x\rangle \xi}\to 0$ in norm. Indeed, note that
$$
\omega_{\delta_n}(a)=(a\langle y_n,x\rangle J\langle y_n,x\rangle \xi,\langle y_n,x\rangle J\langle y_n,x\rangle \xi)
$$
$$
=(a\langle y_n,x\rangle J\langle x,y_n\rangle \langle y_n,x\rangle J\xi,\langle y_n,x\rangle \xi),
$$
and therefore
$$
|(a\langle y_n,x\rangle \xi,\langle y_n,x\rangle \xi)-(a\delta_n,\delta_n)|=|(a\langle y_n,x\rangle (\xi-J\langle x,y_n\rangle \langle y_n,x\rangle J\xi),\langle y_n,x\rangle \xi)| 
$$
$$
\le \|a\| \|\xi-J\langle x,y_n\rangle \langle y_n,x\rangle \xi\|,
$$
which tends to zero. Combining these results one obtains that $\|\omega_{\delta_n}-\omega_{\eta_n}\|\to 0$. Now, because the vectors $\delta_n,\eta_n$ lie in ${\cal P}$, and the fact that norm convergence of vector states with symbols in ${\cal P}$ implies norm convergence of those symbols, one has that $\|\delta_n -\eta_n\|\to 0$ in $H$. In other words,
$$
Re(\langle y_n,x\rangle \xi,J\langle x,y_n\rangle J\eta_n)\to 1.
$$
Suppose now that the states $\omega_{y_n\otimes\eta_n}$ do not converge to $\omega_{x\otimes\xi}$ in the quotient topology of $\Omega_{{\cal A}(X)}$. This means that the fibres of these states do not near in $H$, i.e., there exists a subsequence $y_{n_k}\otimes \eta_{n_k}$ such that 
$\|x\otimes \xi - y_{n_k}\otimes Ju^*\eta_{n_k}\|\ge d > 0$ for all $u\in U_\b$. Or equivalently,
$$
 Re(\langle y_{n_k},x\rangle \xi,Ju^*\eta_{n_k})\le 1-d^2/2, \ \hbox{for all } u\in U_\b .
$$
Clearly this inequality is preserved by taking convex combinations of unitaries $u\in U_\b$ (and leaving everything else fixed), as well as by taking norm limits of such combinations. It follows, using the Russo-Dye theorem, that for $a\in \b$, $\|a\|\le 1$,
$$
Re(\langle y_{n_k},x\rangle \xi,JaJ\eta_{n_k})\le 1-d^2/2.
$$
This clearly contradicts the inequality above, taking $a=\langle x,y_{n_k}\rangle $ for appropriate $k$.
\end{proof}

The next result uses part of the proof of lemma 3. of \cite{pota}
\begin{theo} \label{corolario}
If $\b$ is a separable factor of type II$_1$ such that the tensor product $\b \otimes B(K)$ ($K$ a separable Hilbert space) admits a one parameter automorphism group $\{\theta_s: s\in \zR\}$ scaling the trace of $\b\otimes B(K)$, i.e. $\tau\circ\theta_s=e^{-s}$, $s\in \zR$, with $\tau$ a faithful semi-finite normal trace in $\b \otimes B(K)$, then the map
$$
\vec{\wp_1}:{\cal A}(X)\to \Omega_{{\cal A}(X)} , \vec{\wp_1}(x\otimes\xi)=\omega_{x\otimes \xi}
$$
admits a (global) continuous cross section when $\Omega_{{\cal A}(X)}$ is endowed with the norm topology.
\end{theo}
\begin{proof}{Proof.}
In this case, since $\b$ is finite, $U_\b$ is complete in the strong (=strong$^*$) operator topology \cite{tak}. Moreover, Popa and Takesaki proved in \cite{pota} that it admits a geodesic structure in the sense of Michael \cite{mic}. It has been already remarked that the set function $\omega_{x\otimes \xi}\mapsto \{xu\otimes u^*Ju^*J\xi:u\in U_\b\}$ is lower semicontinuous in the norm topology.  Therefore theorem 5.4 of \cite{mic} applies, and $\vec{\wp_1}$ has a continuous  cross section.
\end{proof}
\begin{coro} 
If $\b$ is a II$_1$ factor satisfying the conditions of \ref{corolario},
then for all  $n\ge 0$, $x\in \s1$, $\f=\omega_\xi \in \Sigma_1(\b)$,
$$
\pi_n(\Omega_{{\cal A}(X)}, \omega_{x\otimes \xi})=\pi_n(\s1,x),
$$
where $\Omega_{{\cal A}(X)}$ is considered with the norm topology, and $\s1$ with the w$^*$ topology.
\end{coro}
\begin{proof}{Proof.}
In \cite{pota} it was proven that the unitary group $U_\b$ of such a factor is contractible
in the ultra strong operator topology, and therefore also in the strong operator topology. The result follows using the above result,
recalling that the fibre of the  fibration $\vec{\wp_1}$ is $U_\b$ with this
topology.
\end{proof}
In \cite{pota} it is noted that remarkable examples of II$_1$ factors
enjoy this property (of having a one parameter group of automorphisms that scale the trace
when tensored with an infinite type I factor), for example ${\cal R}_0$ the hyperfinite II$_1$ factor.

\section{States of the hyperfinite II$_1$ factor}

We will apply the results of the previous section to obtain our main result, namely, that the set of states of ${\cal R}_0$, or more generally, of a factor satisfying the hypothesis of \ref{corolario}, having support equivalent to a given projection $p$, considered with the norm topology, has trivial homotopy groups of all orders.

There is a first result which can be obtained directly from the previous section. 
If ${\cal R}$ is a factor as in \ref{corolario}, and $p\in {\cal R}$ is a proper projection, put $X={\cal R}p$ and $\b=p{\cal R}p$. Clearly $\b$ is a factor which also verifies the hypothesis \ref{corolario}.
Note that $\langle X,X \rangle = span\{px^*yp:x,y\in {\cal R}p\}=p{\cal R}p=\b$ in this case. Therefore by 2.2 of \cite{pas2}, $\{\theta_{x,y}: x,y \in X\}$ spans an ultraweakly dense two sided ideal of $\l$. On the other hand, it is clear that ${\cal R} \subset \l$ as left multipliers, and also that $\theta_{x,y} \in {\cal R}$, for $x,y \in X={\cal R}p$. Indeed, $\theta_{x,y}(z)=x\langle y,z\rangle =xpy^*z$, i.e. left multiplication by $xpy^* \in {\cal R}$. Therefore $\l={\cal R}$. In particular, if $x\in \s1$, $e_x=\theta_{x,x}=xpx^*$ which is equivalent to $px^*xp=\langle x,x\rangle=p$ in $\l$.
The set $\Sigma_{1,X}=\Omega_{{\cal A}(X)}$ equals then the set of states of ${\cal R}$ with support (unitarily) equivalent to $p$. Note that this set is (arcwise) connected in the norm topology. Indeed, if $\b$ is finite, $\s1$ is connected. Using the map $\wp_1$ of section 4, it follows that any two points in $\Sigma_{1,X}$ can be joined with a path (in $\Sigma_{1,X}$) continuous in the $d$-topology, and therefore also in the norm topology.

Applying \ref{corolario} in this situation implies the following:
\begin{coro}
Let ${\cal R}$ be a factor as in \ref{corolario}, and $p\in {\cal R}$ an arbitrary projection. The set of states of ${\cal R}$ with support equivalent to $p$ considered with the norm topology has the same homotopy groups as the set 
$$
{\cal S}_p({\cal R})=\{v\in {\cal R}: v^*v=p\}\subset {\cal R}
$$
regarded with the (relative) ultraweak topology.
\end{coro}
\begin{proof}{Proof.}
In this case $\s1$ clearly equals ${\cal S}_p({\cal R})$ above, and the topology is the w$^*$ (i.e.) ultraweak topology of ${\cal R}$. If $p=0$ the statement is trivial. If $p=1$ it follows from the strong operator contractibility of $U_{\cal R}$ for such ${\cal R}$ proved in \cite{pota}. The case of a proper projection follows from \ref{corolario} and the above remark.
\end{proof}

If $p=0,1$, then ${\cal S}_p({\cal R})$ is contractible (if $p=1$, ${\cal S}_p({\cal R})=U_{\cal R}$). A natural question would be if ${\cal S}_p({\cal R})$ is contractible for proper $p\in {\cal R}$.

We need the following elementary fact:

\begin{lemm} 
Let ${\cal M}\subset B(H)$ be a finite von Neumann algebra, and let $a_n \in {\cal M}$ such that $\|a_n\|\le 1$ and  $a_n^*a_n$ tends to $1$ in the strong operator topology. Then there exist unitaries $u_n$ in ${\cal M}$ such that $u_n-a_n$ converges strongly to zero.
\end{lemm}
\begin{proof}{Proof.}
Consider the polar decomposition $a_n=u_n |a_n|$, where $u_n$ can be chosen unitaries because ${\cal M}$ is finite. Note that $|a_n|\to 1$ strongly. Indeed, since $\|a_n\|\le 1$, $a_n^*a_n\le (a_n^*a_n)^{1/2}$. Therefore, for any unit vector $\xi \in H$, $1\ge (|a_n|\xi,\xi)\ge(a_n^*a_n\xi,\xi)\to 1$. Therefore
$$
\|(a_n-u_n)\xi\|^2=\|u_n(|a_n|-1)\xi\|^2\le \||a_n|\xi- \xi\|^2=1+(a_n^*a_n\xi,\xi)-2(|a_n|\xi,\xi),
$$
which tends to zero.
\end{proof}
In \cite{acs2} it was proven that for a fixed $x_0\in \s1$ the map $\pi_{x_0}:U_{\l} \to \s1$ given by $\pi_{x_0}(U)=U(x_0)$ is onto when $\b$ is finite. In that paper it was considered with the norm topologies. Here we shall regard it with the weak topologies and in the particular case at hand, namely $X={\cal R}p$ and $\b=p{\cal R}p$ with ${\cal R}$ as above. Then, choosing $x_0=p \in \s1={\cal S}_p({\cal R})$, the mapping $\pi_p$ is
$$
\pi_p :U_{\cal R} \to {\cal S}_p({\cal R}),  \ \pi_p(u)=up.
$$
\begin{theo}
If ${\cal R}$ is a factor satisfying the hypothesis of \ref{corolario}, then the map $\pi_p$ above is a trivial principal bundle, when $U_{\cal R}$ is regarded with the strong operator topology and ${\cal S}_p({\cal R})$ is regarded with the ultraweak topology. The fibre is (homeomorphic to) the unitary group of $q{\cal R}q$, where $q=1-p$, again with the strong operator topology.
\end{theo}
\begin{proof}{Proof.}
The key of the argument is again Lemma 3 of \cite{pota}. In that result it is shown that the homogeneous space $U_{\cal R}/U_{\cal M}$ admits a global continuous cross section,  where ${\cal M}\subset {\cal R}$ are factors satisfying the hypothesis of \ref{corolario}, and their unitary groups are endowed with the strong operator topology. In our situation, the fibre of $\pi_p$ (over $p$) is the set $\{u\in U_{\cal R}:up=p\}=\{qwq +p: qwq \in U_{q{\cal R}q}\}= U_{q{\cal R}q}\times \{p\}$. The fibre is not the unitary group of a subfactor with the same unit, nevertheless the argument carries on anyway.
Therefore in order to prove our result it suffices to show that in ${\cal S}_p({\cal R})$  the ultraweak topology (equal to the weak operator topology) coincides with the quotient topology induced by the map $\pi_p$. In other words, that the bijection
$$
U_{\cal R} /U_{q{\cal R}q}\times \{p\} \to {\cal S}_p({\cal R}), \  [u] \to up
$$
is a homeomorphism in the mentioned topologies. It is clearly continuous. It suffices to check continuity of the inverse at the point $p$. Suppose that $u_\alpha$ is a net of unitaries in $U_{\cal R}$ such that $u_\alpha p$ converges weakly to $p$. Then we claim that there are unitaries $qw_\alpha q$ in $q{\cal R}q$ such that $qw_\alpha q+p-u_\alpha$ converges strongly to zero, which would end the proof. This amounts to saying that there exists  unitaries $qw_\alpha q$ verifying that 
$$
Re((qw_\alpha q+p)\xi, u_\alpha \xi)\to \|\xi\|^2
$$
for all $\xi \in H$.
Now since $u_\alpha p \to p$, one has $u_\alpha p\xi \to p\xi$, the former limit is equivalent to the following
$$
Re(qw_\alpha q\xi, u_\alpha q\xi)\to \|q\xi\|^2.
$$
Again, $u_\alpha p\to p$ strongly (and the fact that ${\cal R}$ is finite), imply that $q u_\alpha p$, $pu_\alpha q$, $qu_\alpha^*p$ and  $pu_\alpha^* q$ all converge to zero strongly. Using that $u_\alpha$ are unitaries, these facts imply that $qu_\alpha^*qu_\alpha q\to q$ strongly. Using the lemma above, for the algebra ${\cal M}=q{\cal R}q$, and $a_\alpha=qu_\alpha q$, it follows that there exist unitaries $qw_\alpha q$ in $q{\cal R}q$ such that $qw_\alpha q-qu_\alpha q$ converges to zero strongly. Since $p u_\alpha q$ also tends to zero, it follows that 
$$
qw_\alpha q -u_\alpha q=qw_\alpha q -qu_\alpha q-pu_\alpha q\to 0
$$
strongly. Clearly this last limit  proves our claim.
\end{proof}
Our main result then follows  easily
\begin{theo}
Let ${\cal R}$ be a factor satisfying the hypothesis of \ref{corolario}, and let $p$ be a projection in ${\cal R}$. Then both ${\cal S}_p({\cal R})$ with the ultraweak topology, and the set of normal states of ${\cal R}$ with support equivalent to $p$ with the norm topology, have trivial homotopy groups of all orders.
\end{theo}
\begin{proof}{Proof.}
By the above theorem, ${\cal S}_p({\cal R})$ has trivial homotopy groups, since it is the base space of a fibration with contractible space and contractible fibre. The same consequence holds for the set of normal states with support equivalent to $p$, using the corollary above.
\end{proof}

\begin{rema} \rm
Consider now the restriction of the fibration ${\cal A}(X) \to\Omega_{{\cal A}(X)}$ to the subset
$\{\omega_{x\otimes \xi_0}: x\in \s1\}\subset \Omega_{{\cal A}(X)}$, for a fixed unit, cyclic and separating vector $\xi_0$ i.e.
$$
\{x\otimes \xi_0:x \in \s1\}\simeq \s1 \to \{\omega_{x\otimes \xi_0}: x\in \s1\}, \ \ x\otimes \xi_0\mapsto \omega_{x\otimes\xi_0},
$$
which is again a fibration with the relative topologies. Note that the latter set is in one to one correspondence with  $\of$ of section 2, where $\f=\omega_{\xi_0}$. Therefore one recovers the map $\sigma: \s1 \to \of$, $\sigma(x)=\f_x=\omega_{x\otimes\xi_0}$ of section 2, now considered with the w$^*$ topology for $\s1$ and the norm topology for $\of$. It follows that this map is a fibration, with fibre equal to $U_{{\cal R}^\f}$ with the strong operator topology. 
\end{rema}
One can consider this fibration $\sigma$ in the particular case $X=\b={\cal R}$, for ${\cal R}$ as above, to obtain the following:
\begin{coro}
Let $\f$ be a faithful normal state of a factor ${\cal R}$ as in \ref{corolario}. Then the map
$$
\sigma :U_{\cal R} \to {\cal U}_\f=\{\f\circ Ad(u): u\in U_{\cal R}\}, \ \ \sigma(u)=\f\circ Ad(u)
$$
is a fibration when the unitary group $U_{\cal R}$ is considered with the strong operator topology and the unitary orbit ${\cal U}_\f$ of $\f$ is considered with the norm topology. The fibre is the unitary group $U_{{\cal R}^\f}$ of the centralizer of $\f$ also with the strong operator topology.
Moreover, for $n\ge 0$ one has
$$
\pi_{n+1}({\cal U}_\f,\f)=\pi_n(U_{{\cal R}^\f},1).
$$
\end{coro}
\begin{proof}{Proof.}
It was noted in section 2 that when $X=\b$ is a finite von Neumann algebra, then  $\s1$ is $U_\b$ and $\of$ is the unitary orbit of $\f$.
$\s1=U_\b$ is endowed with the ultraweak topology, which coincides in $U_\b$
with the strong operator topology.
The rest of the corollary follows using that in this case $\sigma$ is (the restriction) of a fibration, and again \cite{pota} that for such factors ${\cal R}$ the unitary group is contractible in the strong operator topology.
\end{proof}

When $n=0$, since $U_{{\cal R}^\f}$ is connected, one obtains that ${\cal U}_\f$ is simply connected in the norm topology. A related result was obtained in \cite{av2}, where it was shown that ${\cal U}_\f$ is simply connected in the quotient topology ($U_\b /U_{\b^\f}$) for any von Neumann algebra $\b$.

Let $p_1,..p_n$ be projections in ${\cal R}$ such that $p_1+...+p_n=1$ and
put $h=r_1p_1+...+r_np_n$ where $r_i$ are positive real numbers such that $r_i\ne r_j$ if $i\ne j$ and $\tau(h)=1$.
Consider the state $\f=\tau (h \cdot  )$. Then clearly ${\cal R}^\f=p_1{\cal R}p_1\oplus ...\oplus p_n{\cal R}p_n$.
Now $U_{p_i{\cal R}p_i}$ is contractible in the strong operator topology,
and therefore $U_{{\cal R}^\f}$ is contractible. It follows that the unitary
orbit ${\cal U}_\f$ (with the norm topology) has trivial homotopy groups of
all orders for such $\f$. Consider this other example: let
${\cal A}\subset {\cal R}$ be a maximal abelian sub (von Neumann) algebra, then
there exists a normal faithful state $\f$ of ${\cal R}$ such that
${\cal R}^\f={\cal A}$. Clearly, since ${\cal R}$ is of type II$_1$, ${\cal A}$
has no atomic projections. It follows that ${\cal A}\simeq L^\infty(0,1)$. It is
fairly elemental to see that $U_{L^\infty(0,1)}$ is contractible in the ultraweak,
i.e. w$^*$ topology. For example consider the map
$F_t(\theta)=\theta \chi_{[0,t)}+\chi_{[t,1]}$, $F_t:U_{L^\infty(0,1)}\to
U_{L^\infty(0,1)}$, for $t \in [0,1]$. Then $F_0(\theta)=\theta$, $F_1(\theta)=1$ for
$\theta \in U_{L^\infty(0,1)}$, $F_t(1)=1$ and $F$ is continuous in the w$^*$ topology:
suppose that $t_\alpha\to t$ and $\theta_\alpha \to \theta$ (w$^*$).
Then if $t_\alpha >t$
$$\int_0^1 [F_{t_\alpha}(\theta_\alpha)(s)-F_t(\theta)(s)]\psi(s)ds=
\int_0^t[\theta_\alpha(s)-\theta(s)]\psi(s)ds+
\int_t^{t_\alpha}[\theta_\alpha(s)-1]\psi(s)ds,
$$
and clearly both integrals tend to zero for any $\psi \in L^1(0,1)$. If $t_\alpha <t$
one reasons analogously. Therefore the map $F_t$ deforms continuously $U_{L^\infty(0,1)}$
to the constant function $1$. It follows that also for such states $\f$, ${\cal U}_\f$ (in the
norm topology) has trivial homotopy groups of all orders. We do not know if
this holds for any faithful normal state of ${\cal R}$.

\vskip2cm

{
%\sit 
\noindent
Esteban Andruchow and Alejandro Varela\\
Instituto de Ciencias \\
Universidad Nacional de Gral. Sarmiento \\
J. A. Roca 850 \\
(1663) San Miguel \\
Argentina  \\
e-mail: eandruch@ungs.edu.ar, avarela@ungs.edu.ar
}

\end{document}